\documentclass[12pt]{article}
\usepackage{enumerate}
\usepackage{amsfonts}
\usepackage{amsmath}
\usepackage{amssymb}

\newcounter{defin}  \newcounter{lemma}  \newcounter{theorem}
\newcounter{property} \newcounter{corol}  \newcounter{remark} \newcounter{example}

\newenvironment{lemma}{\par\refstepcounter{lemma}
     \textbf{Lemma \thelemma.} }{\rm\par}
\newenvironment{theorem}{\par\refstepcounter{theorem}
     \textbf{Theorem \thetheorem.}\ }{\rm\par}
\newenvironment{property}{\par\refstepcounter{property}
     \textbf{Proposition \theproperty.}\ }{\rm\par}
\newenvironment{corollary}{\par\refstepcounter{corol}
     \textbf{Corollary \thecorol.} }{\rm\par}
\newenvironment{definition}{\par\refstepcounter{defin}
     \textbf{Definition \thedefin.}\ }{\rm\par}
\newenvironment{remark}{\par\refstepcounter{remark}
     \textbf{Remark \theremark.}}{\rm\par}
\newenvironment{example}{\par\refstepcounter{example}
     \textbf{Example \theexample.}}{\rm\par}

\begin{document}
\title{Generalized compactness in linear spaces \\
and its applications}

\author{V.Yu.~Protasov\thanks{Lomonosov Moscow State University, e-mail:v-protassov@yandex.ru},
M.E.~Shirokov\thanks{Steklov Mathematical Institute,
e-mail:msh@mi.ras.ru}}

\date{} \maketitle

\vspace{15pt}

\tableofcontents

\section{Introduction}

Various properties and structure of compact sets in the convex
analysis context have been  studied thoroughly starting from the
middle of the last century. An extensive bibliography is devoted to
this theme (see \cite{Phelps, Alf, J&T} and references therein). The
most important results are well known: the Krein-Milman theorem on
convex hulls of extreme points, the Choquet theory on barycentric
decompositions, properties of convex hulls (envelopes) of functions
on convex compact sets. Some of  the classical results have been
extended  to non-compact sets in locally convex spaces by
Edgar~\cite{Edgar, Edgar-2} and  Bourgin~\cite{Bourgin, B&E}. Such
generalizations are interesting not only theoretically, but also
very important in applications, for instance, in mathematical
physics~\cite{M}, in quantum information theory~\cite{Sh-8}, and so
on. Of course, classical results of convex analysis cannot  be
extended to all non-compact sets. One has to postulate some special
properties of these sets. The  Choquet theory, for example, has been
generalized to sets possessing the Radon-Nicodim
property~\cite{Edgar}. In~\cite{Sh-8} the several results on
continuity of convex hulls of functions were extended  to a special
class of sets called $\mu$\nobreakdash-\hspace{0pt}\emph{compact
sets}. This class, characterized by the special relation between the
topology and the structure of linear operations, is the main subject
of this paper.

Problems of continuity of convex hulls of continuous functions (see
definitions in the next section) have been studied in the literature
since the 70s of the last century. Under what conditions on the
convex compact set  $\mathcal A$ is the convex hull of any
continuous (another assumption: concave continuous) function defined
on $\mathcal A$ continuous? Vesterstrom in \cite{Ves} showed that a
necessary and sufficient condition for this is the openness of the
barycenter map. He conjectured the equivalence of the continuity of
the convex hull of any continuous concave function (this property
was called by Lima the
$\mathrm{CE}$\nobreakdash-\hspace{0pt}property~\cite{Lima}) to the
continuity of the convex hull of any continuous function (called
in~\cite{Sh-8} the strong
$\mathrm{CE}$\nobreakdash-\hspace{0pt}property). This conjecture was
proved by O'Brien in \cite{Brien}, who, moreover, showed the
equivalence of the both
$\mathrm{CE}$\nobreakdash-\hspace{0pt}properties to the openness  of
the convex mixing map $$(x, y, \lambda)\mapsto\lambda
x+(1-\lambda)y$$ (the so-called ``stability  property for  convex
sets''~\cite{Susanne, Clausing, Grzaslewicz}). The question arises
if these results can be extended to non-compact sets~$\mathcal A$.
The first step towards the solution of this problem was made
in~\cite{Sh-8}, where so-called
$\mu$\nobreakdash-\hspace{0pt}compact sets were defined. Some
results on the $\mathrm{CE}$\nobreakdash-\hspace{0pt}properties were
generalized from compact sets to the class of
$\mu$\nobreakdash-\hspace{0pt}compact sets. This, in particular,
made it possible to derive several results concerning the entropic
characteristics of infinite-dimensional quantum channels and
systems.

In this paper we analyze the  $\mu$-compactness property in detail,
consider several examples that are important in applications, extend
some classical results of convex analysis known earlier for compact
sets only, in particular, the Vesterstrom-O'Brien theorem, to the
class of~$\mu$-compact sets.

The class of  $\mu$-compact convex  sets is defined by the
requirement of weak compactness of  preimages of all compact sets
under the barycenter map. This property is not purely topological,
it expresses certain relations between topology and the structure of
linear operations. This class contains all compact sets as well as
some important non-compact sets, for example, the set of density
operators in a separable Hilbert space.  $\mu$-compact sets do not
possess many properties of compact sets, such as the boundedness of
continuous functions, the Weierstrass theorem, and so on.
Nevertheless, as we shall see, a lot of results of the Choquet
theory and of the Vesterstrom-O'Brien theory can be extended to this
class. Moreover, we give arguments showing that the class of
$\mu$-compact sets is, in some sense, the largest class, to which
the Vesterstrom-O'Brien theory can be extended.

This paper is organized as follows. In Section~2 we derive basic
properties of $\mu$\nobreakdash-\hspace{0pt}compact sets. By a
simple example we show that several results true for $\mu$-compact
sets become false after a slight relaxing this assumption to {\em
pointwise $\mu$-compactness} (this property is defined by the
requirement of weak compactness of sets of measures with fixed
barycenter). Further we consider examples of $\mu$-compact sets. We
show, in particular, that the bounded part of the positive cone in
the space~~$l_{p}$ for $p=1$ is $\mu$-compact, while for $p>1$ it is
not even pointwise $\mu$-compact. The $\mu$-compactness of the set
of Borel probability measures on a complete separable metric space
is established. This result makes it possible to show that the
convex closure operation respects the $\mu$-compactness property.

In Section~3 we complete the generalization of the
Vesterstrom-O'Brien theorem to the class of $\mu$-compact sets,
started in~\cite{Sh-8}. The $\mu$-compact version of the main result
from~\cite{Brien} is proved. This establishes the equivalence of the
continuity property for convex hulls of \textit{concave} bounded
continuous functions and the continuity property for convex hulls of
arbitrary bounded continuous functions.  We construct an example
confirming our assumption that
 $\mu$\nobreakdash-\hspace{0pt}compact sets form the largest class of convex metrizable sets,
 for which this extension is possible. In Section~4 we apply some of our results to the
 quantum information theory. In Section~5 we discuss possible generalizations and
 formulate several open problems.

\section{On $\mu$-compact sets}

\subsection{Definitions and basic properties}

In Sections 2 and 3 we assume the set  $\mathcal{A}$ to be a closed
bounded subset of a locally convex space. We also suppose that the
convex closure $\overline{\mathrm{co}}\,\mathcal{A}$ of
$\mathcal{A}$ defined as the closure of its convex hull
$\mathrm{co}\,\mathcal{A}$ is a complete separable metric space.
\footnote{This means that the topology on the set
$\overline{\mathrm{co}}\,\mathcal{A}$ is defined by a countable
subset of the family of seminorms generating the topology of the
entire locally convex space, and  this set is separable and complete
in the metric generated by this subset of seminorms.} We use the
following notation:

\begin{description}
  \item[$\mathrm{extr}\, \mathcal{A}$] is the set of extreme
  points of
  $\mathcal{A}$;
  \item[$C(\mathcal{A})$] is the set of continuous bounded functions on a set $\mathcal{A}$;
  \item[$P(\mathcal{A})$ \textmd{and} $Q(\mathcal{A})$] are the sets of convex
  and, respectively,
  concave continuous bounded functions on a convex set  $\mathcal{A}$;
\item[$\mathrm{co}f$ \textmd{and} $\overline{\mathrm{co}}f$] are the convex hull and the convex closure
of a function~$f$ on a convex set; they are defined as the maximal
convex and the maximal convex closed (that is, lower semicontinuous)
functions not exceeding  $f$, respectively \cite{J&T},\cite{P&B};
\item[$\mathfrak{P}_{n}\ =\ \Bigl\{\{\pi_{i}\}_{i=1}^{n}\ |\ \pi_{i}\geq
0,\ \sum_{i=1}^{n}\pi_{i}\ =\ 1\Bigr\}$] is the simplex of all
probability distributions with  $\, n\, \leq \, +\infty$ outcomes.
\end{description}

Let  $M(\mathcal{A})$ be a set of all Borel probability measures on
the set~$\mathcal{A}$ with the topology of weak
convergence~\cite{Bogachev,Par}.


With an arbitrary measure $\mu\in M(\mathcal{A})$ we associate its
barycenter (average)~$\mathbf b\, (\mu)\, \in \,
\overline{\mathrm{co}}\, \mathcal{A}$, which is defined by the
Pettis integral~(\cite{Bogachev, V&T})
\begin{equation}\label{b-exp}
\mathbf b(\mu)\ =\ \int_{\mathcal{A}}x \mu(dx).
\end{equation}

Let  $M_{x}(\mathcal{A})$ be a convex closed subset of the set
$M(\mathcal{A})$, which consists of measures  $\mu$ such that  $\,
\mathbf b(\mu)\, =\, x\, \in\, \overline{\mathrm{co}}\,
\mathcal{A}\, $.

We denote by $\{\pi_{i}, x_{i}\}$ the measure with a finite or
countable number of atoms $\{x_{i}\}$ with weights $\{\pi_{i}\}$.
Let $M^{\mathrm{f}}(\mathcal{A})$ and
$M^{\mathrm{f}}_{x}(\mathcal{A})$ be subsets of the sets
$M(\mathcal{A})$ and $M_{x}(\mathcal{A})$ respectively that consist
of measures with finite supports.

The barycenter map
\begin{equation}\label{b-map}
M(\mathcal{A})\ \ni \ \mu\quad \mapsto \quad \mathbf b(\mu)\ \in \
\overline{\mathrm{co}}\, \mathcal{A},
\end{equation}
is continuous, which can be shown easily by applying Prokhorov's
theorem. Therefore, the image of any compact set in $M(\mathcal{A})$
under map~(\ref{b-map}) is compact in
 $\, \overline{\mathrm{co}}\, \mathcal{A}$.
The inverse map  $\mathbf b^{-1}$ may not possess this property.
Generalizing the definition from~\cite{Sh-8}  consider the class of
convex sets, for which the map  $\mathbf b^{-1}$ takes compact sets
to compact sets.
\smallskip

\begin{definition}\label{ce-def-1}
\textit{The set $\mathcal{A}$ is called $\mu$-compact if the
preimage of any compact subset of $\overline{\mathrm{co}}\,
\mathcal{A}$ under barycenter map~(\ref{b-map}) is a compact subset
of the set~$M(\mathcal{A})$.}
\end{definition}
\smallskip

An arbitrary compact set is $\mu$-compact. Indeed, the compactness
of~$\mathcal{A}$ implies the compactness of~$M(\mathcal{A})$
\cite{Par}. By using Prokhorov's theorem one can derive the
following criterion of $\mu$\nobreakdash-\hspace{0pt}compactness
\cite{Sh-8}.
\medskip

\begin{property}\label{ce-prop-0}
\textit{A convex set $\mathcal{A}$ is $\mu$-compact if and only if
for any compact set $\mathcal{K}\subseteq\mathcal{A}$ and for
any~$\varepsilon>0$ there is a compact subset
$\mathcal{K}_{\varepsilon}\subseteq\mathcal{A}$ such that for any
 $x\in\mathcal{K}$ and for any its expansion  $\, x\, =\,
\sum_{\, i=1}^{\, n}\lambda_{\, i}\, x_{\, i}$, where $\{x_{i}\}_{\,
i=1}^{\, n}\subset\mathcal{A}$, $\{\lambda_{\, i}\}_{\, i=1}^{\,
n}\in\mathfrak{P}_{\, n}$, we have  ${\sum_{\, i\, :\, x_{\, i}\,
\in\, \mathcal{A}\, \setminus\, \mathcal{K}_{\varepsilon}}\,
\lambda_{i}\, <\, \varepsilon}$.}
\end{property}
\medskip

\noindent Proposition~\ref{ce-prop-0} and basic properties of the
set~$M(\mathcal{A})$ yield the following criterion of
$\mu$\nobreakdash-\hspace{0pt}compactness, which is most convenient
for applications.
\smallskip

\begin{property}\label{ce-prop-1}
\textit{A convex set $\mathcal{A}$ is $\mu$-compact if and only if
there is a family~$F(\mathcal{A})$ of nonnegative concave functions
on~$\mathcal{A}$ with the following properties:}
\begin{itemize}
  \item \textit{the set $\{x\in\mathcal{A}\,|\,f(x)\leq c\}$ is
  relatively compact for any function
   $f\in F(\mathcal{A})$ and for any $c>0$;}
  \item \textit{for any compact set
  $\mathcal{K}\subseteq\mathcal{A}$ there is a function
  $f\in F(\mathcal{A})$ such that
  $\sup_{x\in\mathcal{K}}f(x)<+\infty$.}
\end{itemize}
\end{property}

\noindent {\tt Proof.} The sufficiency easily follows from
Proposition~\ref{ce-prop-0} (see \cite{Sh-8}). Let us prove the
necessity. Let $\mathfrak{V}(\mathcal{A})$ be the set of lower
semicontinuous functions~$\varphi$ on $\mathcal{A}$ taking values in
$[0,+\infty]$ and such that  $\{x\in \mathcal{A}\, |\,\varphi(x)\leq
c\}$ is compact for any  $c\geq 0$. Applying Prokhorov's theorem
(see~\cite{Bogachev}, example~8.6.5, p. 236) we conclude that the
set~$M_{0}\subseteq M(\mathcal{A})$ is relatively compact if and
only if there exists a function
$\varphi\in\mathfrak{V}(\mathcal{A})$ such that
$$
\sup_{\mu\in M_{0}}\int_{\mathcal{A}}\, \varphi(x)\, \mu(dx)\ <\
+\infty \, .
$$

Consider the following family of concave nonnegative functions
$$
f_{\, \varphi}(x)=\,\sup_{\mu\in
M_{x}(\mathcal{A})}\int_{\mathcal{A}}\varphi(y)\,
\mu(dy),\;\varphi\, \in\, \mathfrak{V}(\mathcal{A}),
$$
on the set $\mathcal{A}$. This family possesses the first
characteristic property of the family $F(\mathcal{A})$, which
follows from the continuity of the barycenter map. The second
property follows from the~$\mu$-compactness of the
set~$\mathcal{A}$.

{\hfill $\Box$}
\medskip

\begin{remark}\label{ce-prop-0+} This is interesting that
for all convex non-compact, but  $\mu$-compact sets considered in
Subsection~2.2 there exist families $F(\mathcal{A})$ that consist of
affine lower semicontinuous functions.
\end{remark}
\medskip

There exists a criterion of $\mu$-compactness of a convex set in
terms of properties of functions defined on this set~\cite{Sh-10}.
More precisely, it is shown that the $\mu$-compactness is equivalent
to the continuity of the operator of convex closure (that is, the
double Fenchel transform) with respect to monotone pointwise
converging sequences on the classes of continuous bounded and of
lower semicontinuous lower bounded functions.

\medskip

Let us note that continuous affine maps do not necessarily respect
the $\mu$\nobreakdash-\hspace{0pt}compactness property.
Nevertheless, we have the following simple consequence of
Propositions~\ref{ce-prop-0} and~\ref{ce-prop-1}.
\smallskip

\begin{property}\label{mu-compact}
\textit{Let $\mathcal{A}$ and $\mathcal{B}$ be convex
sets,\footnote{It is assumed that the set~$\mathcal{B}$ possess all
the properties mentioned in the beginning of Subsection~2.1.} and
$\varphi$ be a continuous affine map  from $\mathcal{A}$ to
$\mathcal{B}$ such that for any compact set  $\mathcal{C}\subseteq
\mathcal{B}$ its  preimage  $\varphi^{-1}(\mathcal{C})$ is compact
in $\mathcal{A}$. Then}

\textup{1)} \textit{the $\mu$-compactness of~$\mathcal{B}$ implies
the $\mu$-compactness of~$\mathcal{A}$;}

\textup{2)} \textit{if $\varphi$ is surjective, then the
$\mu$-compactness of~$\mathcal{A}$ implies the  $\mu$-compactness
of~$\mathcal{B}$.}
\end{property}
\medskip

The operations of intersection, taking convex closure and Cartesian
product respect $\mu$\nobreakdash-\hspace{0pt}compactness:

\medskip

\begin{property}\label{conv-cl-mu-compact}
1) \textit{A closed subset of any  $\mu$-compact set is
$\mu$-compact.}

2) \textit{The convex closure of any $\mu$-compact set is
$\mu$-compact.}

3) \textit{The Cartesian product of a finite or countable family of
$\mu$-compact sets is $\mu$\nobreakdash-\hspace{0pt}compact (in the
topology of coordinate-wise convergence).}
\end{property}
\smallskip

\noindent {\tt Proof.} 1) This follows directly from
Definition~\ref{ce-def-1}.

2) Combining the $\mu$-compactness of the set $\mathcal{A}$ and
proposition~2 from~\cite{Sh-8} (its proof for our class of sets is
literally the same as in that paper) we obtain that the barycenter
map $\mu\mapsto \mathbf b(\mu)$ is a continuous affine surjection
from $\, M(\mathcal{A})\, $ to $\, \overline{\mathrm{co}}\,
\mathcal{A}$ satisfying the assumptions of
Proposition~\ref{mu-compact}. Applying the second part of that
proposition and Corollary~\ref{p-measures} from the next
subsection we conclude that the set $\, \overline{\mathrm{co}}\,
\mathcal{A}\, $ is $\mu$-compact.

3) By assertion 2) it suffices to consider the case of convex
 $\mu$-compact sets. Assume for every  $\, n \in \mathbb{N}\, $ the set $\, \mathcal A^{\,n}\, $
 is $\mu$-compact. Let us show that the set $\mathcal A \,= \,
\otimes_{n \in \mathbb{N}} \, \mathcal A^{\, n}\,$ is $\mu$-compact
in the topology of coordinatewise convergence. For an  arbitrary
compact set $\, \mathcal K \, \subset \, \mathcal A \, $ and for
each  $\,n \in \mathbb{N}\,$ let  $\, \mathcal{K}^{\, n}\, $ be the
projection of $\, \mathcal K\, $ onto $\, \mathcal A^{\, n}$. The
set $\, \mathcal{K}^{\, n}\, $ consists of points  $\, x^{\,n} \,
\in \, \mathcal A^{\, n}$ which are the corresponding coordinates of
some point $x \in \mathcal K$. This set is compact. Since the set
$\mathcal A^{\, n}\, $ is $\, \mu$-compact, it follows from
Proposition~\ref{ce-prop-0} that for any $\varepsilon
> 0$ there exists the corresponding compact
$\, \, \mathcal K^{\, n}_{\, \varepsilon}\, \subset \, \, \mathcal
A^{\, n}$. Since  $\, \mathcal K \, \subseteq \, \otimes_{n \in
\mathbb{N}} \, \mathcal{K}^{\, n}$, we have  $\, \mathcal K_{\,
\varepsilon}\, = \, \otimes_{n \in \mathbb{N}} \, \mathcal K^{\,
n}_{\varepsilon2^{-n}}$. It is easy to check that this set satisfies
the assumptions of Proposition~\ref{ce-prop-0}. Therefore, the
set~$\mathcal A\, $ is $\, \mu$-compact. $\square$

\medskip

The first assertion of Proposition~\ref{conv-cl-mu-compact} implies
that the intersection of  $\mu$\nobreakdash-\hspace{0pt}compact sets
is $\mu$\nobreakdash-\hspace{0pt}compact. However,  their union or
Minkowski sum is, in general, not $\mu$-compact (Remark~\ref{Mink}).

\begin{remark}\label{conv-cl-mu-compact+}
By the second assertion of Proposition~\ref{conv-cl-mu-compact} to
prove  the $\mu$-compactness of a convex set $\mathcal{A}$ it
suffices to show the $\mu$-compactness of any its subset
$\mathcal{B}$ such that $\, \mathcal{A}\, =\,
\overline{\mathrm{co}}\, \mathcal{B}$. Note that the sets of
measures $M(\mathcal{A})$ and $M(\mathcal{B})$ (that are involved in
the definition of $\mu$-compactness of the sets $\mathcal{A}$ and
$\mathcal{B}$) can be completely different. For example, if
$\mathcal{A}$ is a simplex, then $\mathcal{B}$ can be a countable
family $\{e_{i}\}$ of isolated extreme points of the
set~$\mathcal{A}$. Hence $M(\mathcal{B})$ is isomorphic to the set
$\mathfrak{P}_{+\infty}$ of all probability distributions with
countable number of outcomes. The criterion of $\mu$-compactness for
the set $\mathcal{B}$, and therefore, for the set~$\mathcal{A}$, can
be formulated as follows: for any compact set
$\mathcal{K}\subset\mathcal{A}$ and for any $\varepsilon>0$ there
exists $n$ such that the inclusion
$\sum_{i=1}^{+\infty}\pi_{i}e_{i}\in\mathcal{K},\;\{\pi_{i}\}\in\mathfrak{P}_{+\infty},$
implies  $\sum_{i=n+1}^{+\infty}\pi_{i}<\varepsilon$.

The second assertion of Proposition~\ref{conv-cl-mu-compact}
together with  Proposition~\ref{KMC} below lead to the following
observation. Let $\mathcal{A}$ be a $\mu$-compact convex set in the
initial topology~$\tau$ and let $\tau'$ be a stronger topology on
$\mathcal{A}$, that coincides with $\tau$ on the set
$\overline{\mathrm{extr}\, \mathcal{A}}$; then the set $\mathcal{A}$
is $\mu$-compact in the topology~$\tau'$.
\end{remark}

Proposition~\ref{mu-compact} implies that all isomorphisms in the
category of $\mu$-compact convex sets are affine homeomorphisms. The
affinity assumption cannot be omitted, which is seen from  the
following example. Suppose $\mathcal{A}$ is a convex set that is not
$\mu$-compact. The first assertion of
Proposition~\ref{conv-cl-mu-compact} and Corollary~\ref{p-measures}
from the next subsection yield that the subset of the set
$M(\mathcal{A})$ that consists of Dirac (single-atomic) measures is
a $\mu$-compact set, which is, moreover, homeomorphic to the
set~$\mathcal{A}$. This observation shows that the $\mu$-compactness
property, in contrast to the compactness, is not purely topological.
It is defined by composition of the topology and of the structure of
the operation of convex mixing.
\medskip

\noindent Proposition~\ref{mu-compact} gives the following condition
of $\mu$-compactness of families of maps. This condition is applied
in the next section.
\medskip

\begin{corollary}\label{map-mu-compact}
\textit{Let $\mathfrak{F}(\mathcal{X},\mathcal{Y})$ be a locally
convex space with the topology~$\tau$ of maps from the set
$\mathcal{X}$ to a locally convex space $\mathcal{Y}$. Let also
$\mathfrak{F}_{0}$ be a convex closed bounded subset of the space
$\mathfrak{F}(\mathcal{X},\mathcal{Y})$ that consists of maps taking
values from the convex $\mu$-compact set
$\mathcal{A}\subset\mathcal{Y}$. Moreover, $\mathfrak{F}_{0}$ is a
complete separable metric space, for which there is an element
$x_{0}\in \mathcal{X}$ such that:}
\begin{enumerate}[1)]
  \item \textit{$\{\,\tau\textrm{-}\lim_{n\rightarrow+\infty}\Phi_{n}=\Phi_{0}\,\}\;\Rightarrow\;
  \{\,\lim_{n\rightarrow+\infty}\Phi_{n}(x_{0})=\Phi_{0}(x_{0})\,\},\quad\forall
  \{\Phi_{n}\}\subset\mathfrak{F}_{0}$;}
  \item
  \textit{the set
  $\,\{\Phi\in\mathfrak{F}_{0}\,|\,\Phi(x_{0})\in\mathcal{C}\}$
  is relatively $\tau$-compact
 for any compact set $\;\mathcal{C}\subseteq\mathcal{A}$.}
\end{enumerate}

\textit{Then the set $\mathfrak{F}_{0}$ is $\mu$-compact.}
\end{corollary}
\vspace{5pt}

\noindent {\tt Proof.} The continuous affine map
$$
\mathfrak{F}_{0}\, \ni \, \Phi\ \mapsto \ \Phi(x_{0})\, \in \,
\mathcal{A}
$$
satisfies the assumptions of Proposition~\ref{mu-compact}. The first
part of that proposition implies the $\mu$-compactness of the set
$\mathfrak{F}_{0}$. $\square$

\medskip

For further analysis of the $\mu$-compactness we need a weaker
version of this property.\vspace{5pt}

\begin{definition}\label{ce-def-2}
\textit{The set $\mathcal{A}$ is called pointwise $\mu$-compact if
for any
 $x\in\overline{\mathrm{co}}\, \mathcal{A}$
the set $M_{x}(\mathcal{A})$ is a compact subset of the set
$M(\mathcal{A})$.}
\end{definition}
\medskip

Clearly, the pointwise  $\mu$-compactness follows from the
$\mu$-compactness. However, as we shall see in
Proposition~\ref{l1-lp}, these two properties are not equivalent.

\bigskip

$\mu$-compact sets do not possess many properties of compact sets,
such as uniform continuity and boundedness of continuous functions,
the Weierstrass theorem of extremal values of continuous functions,
and so on. It appears, however, that $\mu$-compact sets inherit some
important properties of compact sets. This allows us to extend many
results of the  Choquet theory and of the Vesterstrom-O'Brien theory
to those sets (see the next section).

If $\mathcal{A}$ is a  convex set, then one can introduce the
following partial order, called the Choquet order, on the set
$M(\mathcal{A})$  \cite{Phelps},\cite{Bourgin}. We suppose that
$\mu\succ\nu$ if and only if
$$
\int_{\mathcal{A}}f(y)\, \mu(dy) \ \geq\ \int_{\mathcal{A}}f(y)\,
\nu(dy)
$$
for any function  $f\in P(\mathcal{A})$.

A measure $\mu\in M(\mathcal{A})$ is called \textit{maximal} if
$\nu\succ\mu$ implies  $\nu=\mu$ for any measure $\nu\in
M(\mathcal{A})$. If  $\mu$ and $\nu$ are measures from
$M(\mathcal{A})$ such that  ${\mu\succ\nu}$, then $\, \mathbf
b(\mu)\, =\, \mathbf b(\nu)$. This follows from the fact that the
set of continuous bounded affine functions on the set $\mathcal{A}$
separates its points. The following fact is a straightforward
consequence of Definition~\ref{ce-def-2}.
\smallskip

\begin{lemma}\label{max}
\textit{Let  $\mathcal{A}$ be a convex pointwise $\mu$-compact set.
Any subset of the set $M(\mathcal{A})$ that is linearly ordered by
the relation $\,\prec\,$ has the least upper bound.}
\end{lemma}
\smallskip

\noindent {\tt Proof.} Any subset of the set $M(\mathcal{A})$
linearly ordered by the relation $\prec$, which is actually a net
$\{\mu_{\lambda}\}_{\lambda\in\Lambda}$, is contained in
$M_{x}(\mathcal{A})$ for some $x\in\mathcal{A}$ and, consequently,
is relatively compact. Whence there is a subnet
$\{\mu_{\lambda_{\pi}}\}_{\pi\in\Pi}$ that converges to some measure
$\mu_{0}\in M_{x}(\mathcal{A})$. One can easily verify that
$\mu_{\lambda}\prec\mu_{0}$ for all $\lambda\in\Lambda$.

{\hfill $\Box$}
\medskip

Combining Lemma~\ref{max} with theorem~2.4 from~\cite{Edgar-2} we
see that the class of convex pointwise  $\mu$-compact sets is a
proper subclass of the class of convex sets with the Radon-Nicodim
property. This class is studied in an extensive
literature~\cite{Edgar,Edgar-2,Bourgin}.

Applying Lemma~\ref{max} and Zorn's lemma one can easily prove the
following generalizations  of the Krein-Milman theorem and of the
Choquet theorem to pointwise $\mu$-compact sets. This actually
follows from the Radon-Nicodim property of these
sets~\cite{Edgar-2}.

\medskip

\begin{property}\label{KMC}
\textit{Let $\mathcal{A}$ be a pointwise $\mu$-compact convex set.
Then $\ \overline{\mathrm{co}}\, (\mathrm{extr}\, \mathcal{A})\ =\
\mathcal{A}\, $ and \linebreak $\, {\mathbf b\left(M\left(\,
\overline{\mathrm{extr}\, \mathcal{A}}\, \right)\right)\ =\
\mathcal{A}}$.}
\end{property}
\smallskip

\noindent {\tt Proof.} Since the set of finitely supported measures
is dense in~$M(\overline{\mathrm{extr}\, \mathcal{A}})$
\cite[theorem 6.3]{Par}, we see that the first assertion follows
from the second one.

Let $x_{0}\in\mathcal{A}$. Lemma~\ref{max} combined with Zorn's
lemma yields the existence of the maximal measure $\mu_{*}$ in
$M_{x_{0}}$. From the properties of the Choquet ordering it follows
that the measure $\mu_{*}$ is also maximal in $M(\mathcal{A})$;
hence it lies in $M\left(\, \overline{\mathrm{extr}\, \mathcal{A}}\,
\right)$ (see the proof of theorem~5.2 in~\cite{Edgar}). Thus, $\,
x_{0}\, \in \, \mathbf b(M\left(\, \overline{\mathrm{extr}\,
\mathcal{A}}\, \right))$.

{\hfill $\Box$}
\medskip


Another property inherited by $\mu$-compact sets from compact ones
is the following representation of convex closure of lower
semicontinuous functions \cite[proposition 3]{Sh-8}. Its proof is
easily generalized to our class of sets.
\medskip

\begin{property}\label{formula}
\textit{If  $f$ is a lower semicontinuous lower bounded function on
a convex $\mu$\nobreakdash-\hspace{0pt}compact set~$\mathcal{A}$,
then its convex closure can be determined by the expression}
\begin{equation}\label{co-cl-exp}
\overline{\mathrm{co}}\, f(x)\ = \ \inf_{\mu\in
M_{x}(\mathcal{A})}\int_{\mathcal{A}} f(y) \, \mu(dy)\, , \qquad
\forall x\, \in \, \mathcal{A},
\end{equation}
\textit{where the infimum is attained at some measure $\mu^{f}_{x}$
from  $M_{x}(\mathcal{A})$.}
\end{property}
\smallskip

This representation is a crucial point for most results on convex
closures of functions defined on~$\mu$-compact sets.
\medskip

If $f$ is a continuous bounded function on a convex~$\mu$-compact
set  $\mathcal{A}$, then combining the continuity of the functional
$M(\mathcal{A})\ni\mu\mapsto\int_{\mathcal{A}}f(x)\mu(dx)$
 and \cite[lemma 1]{Sh-8} we see that the infimum
in~(\ref{co-cl-exp}) can be taken over finitely-supported
measures.\medskip

\begin{corollary}\label{formula+}
\textit{An arbitrary continuous bounded function $f$ on a convex
$\mu$-compact set $\mathcal{A}$ possesses a lower semicontinuous
(closed) convex hull, that is, }
\begin{equation}\label{co-cl-exp+}
\mathrm{co}f(x)=\inf_{\{\pi_{i},x_{{i}}\}\in
M^{\mathrm{f}}_{x}(\mathcal{A})}\sum_{i}\pi_{i}f(x_{i})\ =\
\overline{\mathrm{co}}f(x),\quad \forall x\in \mathcal{A}.
\end{equation}
\end{corollary}
This is a  $\mu$-compact generalization of corollary~I.3.6 in
\cite{Alf}. The $\mu$-compactness assumption about the set
$\mathcal{A}$ in Proposition~\ref{formula} and in
Corollary~\ref{formula+} cannot be weakened to pointwise
$\mu$\nobreakdash-\hspace{0pt}compactness.\medskip

\begin{property}\label{formula-}
\textit{If a convex set $\mathcal{A}$ is not $\mu$-compact, but
pointwise  $\mu$-compact, then there exists a continuous bounded
function $f$ on $\mathcal{A}$, whose convex hull $\,\mathrm{co}f$ is
not lower semicontinuous.} \textit{This means that
representation~(\ref{co-cl-exp}) for the convex closure
$\,\overline{\mathrm{co}}f$ of the function~$f$ does not hold.}
\end{property}\medskip

\noindent {\tt Proof.} Since the set $\mathcal{A}$ is not
$\mu$-compact, there exists a non-compact sequence of measures
$\{\mu_{n}\}\subset M(\mathcal{A})$ such that the corresponding
sequence $\{x_{n}=\mathbf{b}(\mu_{n})\}\subset\mathcal{A}$
converges. By Prokhorov's theorem the sequence $\{\mu_{n}\}$ is not
tight. As it was shown in the proof of theorem~8.6.2 in
\cite{Bogachev} this guarantees the existence of $\varepsilon>0$ and
$\delta>0$ such that for any compact set
$\mathcal{K}\subset\mathcal{A}$ and any natural $N$ there is $n>N$,
for which $\mu_{n}(U_{\delta}(\mathcal{K}))<1-\varepsilon$, where
$U_{\delta}(\mathcal{K})$ is a closed $\delta$-neighborhood of the
compact set $\mathcal{K}$. Since finitely supported measures is
dense in the set of all measures with a fixed barycenter~\cite[lemma
1]{Sh-8}, applying corollary~8.2.9~from~\cite{Bogachev}, we may
assume that the sequence $\{\mu_{n}\}$ contains only finitely
supported measures.

Let $x_{0}$ be a limit of the sequence~$\{x_{n}\}$. By the pointwise
$\mu$-compactness of the set $\mathcal{A}$ for the value
$\varepsilon$ defined above there exists a convex compact set
$\mathcal{K}_{\varepsilon}$ such that
$\mu(\mathcal{K}_{\varepsilon})\geq 1-\varepsilon/2$ for any measure
$\mu\in M_{x_{0}}(\mathcal{A})$.

Let $f$ be a bounded continuous function on  $\mathcal{A}$ such that
$f(x)=1$ for all $x\in\mathcal{K}_{\varepsilon}$ and $f(x)\leq 0$
for every $x\in\mathcal{A}\setminus
U_{\delta}(\mathcal{K}_{\varepsilon})$. Then from the properties of
the sequences $\{x_{n}\}$ and $\{\mu_{n}\}$ it follows that for any
natural $N$ there is $n>N$, for which
$\mathrm{co}f(x_{n})<1-\varepsilon$. On the other hand, the
properties of the set $M_{x_{0}}(\mathcal{A})$ imply that
$\mathrm{co}{f}(x_{0})\geq 1-\varepsilon/2$. Thus, the function
$\mathrm{co}f$ is not lower semicontinuous.

{\hfill $\Box$}
\medskip

\begin{remark}\label{r3}
The pointwise $\mu$-compactness condition cannot be omitted in the
proof of Proposition~\ref{formula-}. This shows that $\mu$-compact
sets do not form the maximal class of convex sets, for which
Corollary~\ref{formula+} holds. However, without the pointwise
$\mu$-compactness Corollary~\ref{formula+} may fail (see
Example~\ref{conv-cl-ex} in the next section).

To construct an example showing that the pointwise~$\mu$-compactness
condition is essential in Proposition~\ref{formula-}, we consider
the unit ball~$\mathcal B$ in a separable Hilbert space. Clearly,
that is not a pointwise $\mu$-compact set. Let us show that {\em for
any continuous bounded function  $f$ on~$\mathcal B$ its convex hull
$\mathrm{co}f$ is continuous.} Since a bounded convex function is
continuous at interior points of its domain, it suffices to prove
the continuity of  $\mathrm{co}f$ on the boundary of $\mathcal B$,
that is, on the unit sphere.

Let $x$ be a point of the unit sphere. For  arbitrary~$\varepsilon
>0$ there exists $\delta >0$ such that $\, |f(x) - f(y)|\, < \, \varepsilon$, whenever $\, \|x - y\| < \, 2\,
\delta\, $. Applying Lemma~\ref{hilbertball} stated below, we
conclude that for any $z$ such that $\, \|z-x\|\, < \, \delta\, $
and for any $\mu \in M_z(\mathcal B)$ the measure of a ball of
radius $2\, \delta$ centered at~$x$ is at least $\, r(\delta , z)$.
Therefore, the value $\, \bigl|\, f(x) \, - \, \int_{\mathcal B}f(t)
\, d \mu\, |\, $ does no exceed $\, \varepsilon \, r(\delta , z)\, +
\, N \, \bigl(\, 1\, - \, r(\delta , z) \, \bigr)$, where $\, N \, =
\, \sup_{t\in \mathcal B}|f(t)|$. Whence
$$
\bigl|\, f(x) \, - \, \mathrm{co}f(z) \, | \quad \le \quad
\varepsilon \, r(\delta , z)\, + \, N \, \bigl(\, 1\, - \, r(\delta
, z) \, \bigr).
$$
Since for $\, z \to x$ we have  $\|z\| \, \to \, 1\, $, it follows
that$\, r(\delta , z) \, \to \, 1\, $ and so $\, \mathrm{co}f(z)\,
\to \, f(x)$. It remains to note that $\mathrm{co}f(x)\,  = \,
f(x)$, because $\, x$ is an extreme point of $\mathcal B$.
\end{remark}

\medskip
\begin{lemma}\label{hilbertball}
\textit{Let  $\mathcal B$ be the unit ball of the Hilbert space and
let~$\, \delta
>0$. For an arbitrary point $z \in \mathcal B$ such that $\, \|z\| \, >
\, 1\, - \, \delta $ and for any measure $\mu \in M_{z}(\mathcal B)$
the measure of a ball of radius $\delta$ centered at $z$ is at least
$\, r(\delta , z) \, = \, \frac{\delta^2 \, - \,
(1-\|z\|^2)}{\delta^2 \, - \, (1-\|z\|)^2}$.}
\end{lemma}\medskip
\noindent {\tt Proof}.  Let $\, \bar z\, = \, \frac{z}{\|z\|}\, $
and
$$
\mathcal B_0 \ = \ \left\{\, y \in \mathcal B,\ (y,\bar z) \, < \,
\frac{1+\|z\|^2 - \delta^2}{2\, \|z\|}  \, \right\},\qquad \mathcal
B_1 \ = \ \left\{\, y \in \mathcal B,\ (y,\bar z) \, \ge
 \, \frac{1+\|z\|^2 - \delta^2}{2\,
\|z\|}  \, \right\}.
$$
By a direct calculation we show that if  $\, \|y\|\, = \, 1\, $ and
$\, \|y-z\| > \delta$, then $\, y \in \mathcal B_0$. Consequently,
all points of the ball $\mathcal B$ lying at a distance more than
$\delta$ from the point  $z$ are located in the set $\mathcal B_0$.
Denote by  $c_i$ the barycenter of the measure $\, \mu$ on the set
$\mathcal B_i\, , \, i = 0,1$. We have $\, z \, = \, \mu(\mathcal
B_0)\, c_0 \, + \, \mu(\mathcal B_1) \, c_1.$ It is obvious that $\,
(c_0, \bar z)\, \leq \, \frac{1+\|z\|^2 - \delta^2}{2\, \|z\|}\, $
and  $\, (c_1, \bar z) \, \le \, 1$. Therefore,
$$
\|z\|\, = \, (z, \, \, \bar z) \ < \ \bigl( 1\, -\, \mu(\mathcal
B_1)\bigr)\, \frac{1+\|z\|^2 - \delta^2}{2\, \|z\|}\ + \
\mu(\mathcal B_1),
$$
which gives the desired inequality for the measure of~$\mathcal
B_1$, and hence for the measure of the ball of radius~$\delta $
centered at~$z$, since this ball contains~$\, \mathcal B_1$.

{\hfill $\Box$}
\medskip

\medskip

One of the most important conclusions from Proposition~\ref{formula}
is that any lower semicontinuous lower bounded function~$f$ defined
on a convex $\mu$-compact set $\mathcal{A}$ coincides with its
convex closure $\overline{\mathrm{co}}f$ on the set of extreme
points $\mathrm{extr}\, \mathcal{A}$. Furthermore, that proposition
enables us to obtain the following representation for the
set~$\mathrm{extr}\, \mathcal{A}$, which will help us in Section~3.
\medskip

\begin{property}\label{ext-rep}
\textit{Let $\mathcal{A}$ be a $\mu$-compact convex set. Then}
\begin{equation}\label{ext-rep+}
\mathrm{extr}\, \mathcal{A}\, =\, \bigcap_{f\in
Q(\mathcal{A})}\mathcal{B}_{f},\quad \textit{where} \quad
\mathcal{B}_{f}\, =\,
\{x\in\mathcal{A}\,|\,f(x)=\overline{\mathrm{co}}f(x)\}.
\end{equation}
\end{property}
\smallskip

\noindent {\tt Proof.} The inclusion $\mathrm{extr}\, \mathcal{A}\,
\subseteq\mathcal{B}_{f}$ for any function~$f\in Q(\mathcal{A})$
follows from~(\ref{co-cl-exp+}). Suppose
$x_{0}\in\mathcal{A}\setminus\mathrm{extr}\, \mathcal{A}$; then
there are two distinct points  $x_{1}$ and $x_{2}$ in $\mathcal{A}$
such that $x_{0}=\frac{1}{2}x_{1}+\frac{1}{2}x_{2}$. To prove that
$x_{0}$ is not in $\bigcap_{f\in Q(\mathcal{A})}\mathcal{B}_{f}$ one
needs to find a function~$f$ from $Q(\mathcal{A})$ such that
$f(x_{0})>\frac{1}{2}f(x_{1})+\frac{1}{2}f(x_{2})$. The function
$-a^{2}(\cdot)$ will do, where $a$ is an affine continuous bounded
function on $\mathcal{A}$ such that $a(x_{1})\neq a(x_{2})$. {\hfill
$\Box$}
\medskip

Example~\ref{conv-cl-ex} in the next subsection shows the importance
of the $\mu$-compactness assumption in Proposition~\ref{ext-rep}.

Arguing as in the proof of Proposition~\ref{formula-} one can easily
produce the following necessary condition for
representation~(\ref{ext-rep+}). This is a condition of local
$\mu$-compactness in neighborhood of the set~$\mathrm{extr}\,
\mathcal{A}$: \textit{for any extreme point~$x_{0}$ of the
set~$\mathcal{A}$ and for any sequence~$\{x_{n}\}\subset\mathcal{A}$
converging to~$x_{0}$ the set $\mathbf{b}^{-1}(\{x_{n}\})$ is
compact in $M(\mathcal{A})$}.

\subsection{Examples}

Any  compact set is obviously $\, \mu$-compact. In this section we
consider the several most important examples of $\, \mu$-compact,
but not compact sets.\medskip

\begin{property}\label{V1}
\textit{The bounded part\footnote{Here and in what follows by the
bounded part of a positive cone in an ordered Banach space we mean
the intersection of this cone with a unit ball.} of the positive
cone in the space
 $l_{1}$ is
$\mu$-compact.}
\end{property}
\smallskip

\noindent {\tt Proof}. It is sufficient to take the family of
functions of the form $l_{1}\ni\{x_{i}\}\mapsto \sum_{i}h_{\, i}\,
x_{\, i}$, where $\{h_{\, i}\}$ is an increasing unbounded sequence
of positive numbers, and then apply Proposition~\ref{ce-prop-1}
taking into account the compactness criterion for subsets of the
space~$l_{1}$.$\;\Box$

\medskip

\begin{corollary}\label{prob}
\textit{The set $\,\mathfrak{P}_{+\infty}$ of all probability
distributions with countably many outcomes is a $\mu$-compact subset
of the space~$\, l_{\, 1}$. }
\end{corollary}
\medskip

\begin{remark}\label{Mink}
Let  $\, \mathcal{A}_1\, = \, \bigl\{\, x \in l_{\,1} \ \bigl| \ x
\ge 0\, , \|x\|_{\, 1}\, \le \, 1\, \bigr\}$. Proposition~\ref{V1}
implies that the sets $\mathcal{A}_1$ and $-\mathcal{A}_1$ are both
$\, \mu$-compact in the metric $\, l_{\, 1}$. However, neither their
convex hull nor their Minkowski sum is $\mu$-compact. They are
actually not even pointwise $\mu$-compact, because they contain a
unit ball of the space~$\,l_{\,1}$, which is not pointwise
$\mu$-compact (this is easy to show).
\end{remark}

\medskip

\begin{property}\label{measures}
\textit{An arbitrary weakly closed  bounded (in variation) set of
Borel measures on a complete separable metric space is $\mu$-compact
in the weak convergence topology.}
\end{property}
\medskip

\noindent {\tt Proof.} Let $X$ be a complete separable metric space,
$M$ be a weakly closed bounded set of Borel measures on that
space~$X$, $\, \mathfrak{V}(X)$ be the set of all lower
semicontinuous functions~$\varphi$ on~$X$ taking values
from~$[0,+\infty]$ such that the set $\{x\in X\, |\,\varphi(x)\leq
c\}$ is compact for all  $c\geq 0$. Prokhorov's theorem  (see
\cite{Bogachev}, example 8.6.5, p. 236) implies that the set
$M_{0}\subseteq M$ is relatively compact if and only if there exists
a function~$\varphi\in\mathfrak{V}(X)$, for which
$$
\sup_{\mu\in M_{0}}\int_{X}\, \varphi(x)\, \mu(dx)\ <\ +\infty\, .
$$
Hence the family of affine lower semicontinuous functions $\,
f_{\varphi}(\mu) = \int_{X}\varphi(x) \mu(dx)$, $\varphi \in
\mathfrak{V}(X)$ on the set~$M$ satisfies all the assumptions of
Proposition~\ref{ce-prop-1}. {\hfill $\Box$}
\medskip

\begin{corollary}\label{p-measures}
\textit{The set of all Borel probability measures on a complete
separable metric space is~$\mu$-compact in the topology of weak
convergence.}
\end{corollary}\medskip

Using Propositions~\ref{mu-compact} and \ref{V1} as well as
Proposition~\ref{lp} stated below one can prove the following
assertion on the properties of the cone of positive operators in the
Shatten class of order $p$, that is, in the Banach space of all
operators acting in a separable Hilbert space $\mathcal{H}$ such
that $\mathrm{Tr}|A|^{p}<+\infty$ with the norm $\, \|A\|_{p}\, =\,
\bigl( {\,\mathrm{Tr}\, |A|^{p}}\ \bigr)^{1/p}\, $.
\medskip

\begin{property}\label{Lp}
\textit{The bounded part of the positive cone in the Shatten class
of order $p$ is $\mu$-compact precisely for $p=1$.}
\end{property}
\smallskip

\noindent {\tt Proof.} In the case $p=1$ by applying the compactness
criterion for subsets of the positive
cone~$\mathfrak{T}_{+}(\mathcal{H})$ of the space of trace class
operators $\mathfrak{T}(\mathcal{H})$ (\cite{Sh-H-3}, Appendix), we
obtain that the map taking an operator
$A\in\mathfrak{T}_{+}(\mathcal{H})$ to the sequence of its diagonal
elements (in some fixed basis of the space $\mathcal{H}$) satisfies
the assumptions of the first part of Proposition~\ref{mu-compact}
(with the bounded parts of the positive cones of the spaces
$\mathfrak{T}(\mathcal{H})$ and $l_{1}$ in the roles of
$\mathcal{A}$ and $\mathcal{B}$ respectively). Combining that
proposition with Proposition~\ref{V1} we obtain the
$\mu$-compactness of the bounded part of the
cone~$\mathfrak{T}_{+}(\mathcal{H})$.

The cone of positive operators in the Shatten class of order $p>1$
contains a subcone of commuting operators. That subcone is affinely
homeomorphic to $\mathcal{A}_p$ (the bounded part of the positive
cone in the space~$l_{p}$), which  is not
 $\mu$-compact  (Proposition~\ref{lp}).$\;\Box$

\medskip

Proposition~\ref{Lp} yields, in particular, the $\mu$-compactness of
the set~$\mathfrak{S}(\mathcal{H})$ of quantum states --  density
operators in a separable Hilbert space $\mathcal{H}$. This fact was
originally proved in~\cite{H-Sh-2}. Corollary~\ref{map-mu-compact}
and Lemma~\ref{comp-c} from Appendix imply the following
result.\medskip

\begin{property}\label{p-map}
\textit{Let $\mathfrak{L}_{+}(\mathcal{H},\mathcal{H}')$ be the cone
of positive linear continuous maps from the Banach
space~$\mathfrak{T}(\mathcal{H})$  of trace class operators in a
separable Hilbert space~$\mathcal{H}$ into a similar
space~$\mathfrak{T}(\mathcal{H'})$. Then the bounded (in the
operator norm) part of this cone is $\mu$-compact in the strong
operator topology.}
\end{property}\medskip

Proposition~\ref{p-map} implies  $\mu$-compactness of the sets of
quantum operations and of quantum channels in the topology of strong
convergence~\cite{Sh-H-3}.

Proposition~\ref{p-map} also yields the $\mu$-compactness (in the
strong operator topology) of the bounded part of the  cone of linear
continuous positive operators in~$l_{1}$. Indeed, this cone is
naturally identified with a subset of the cone
$\mathfrak{L}_{+}(\mathcal{H},\mathcal{H}')$. This shows the
$\mu$-compactness of the set of sub-Markov operators, in particular,
Markov operators.
\medskip

Now consider some ``negative'' examples. Proposition~\ref{l1-lp}
gives an example of pointwise $\mu$-compact sets
(Definition~\ref{ce-def-2}) that are not $\mu$-compact. We are going
to see that in Hilbert space there are no $\mu$-compact sets that
are not compact (Proposition~\ref{l2}). Then we give several
examples of sets that are not even pointwise $\mu$-compact.
\medskip

By definition the pointwise $\mu$-compactness property of a set
survives weakening the topology. The next proposition shows that
$\mu$-compactness does not possess this property. The bounded part
of the positive cone in~$l_{1}$ loses the $\mu$-compactness after
weakening the topology. \smallskip
\begin{property}\label{l1-lp}
\textit{For any $p>1$ the simplex $\,\Delta_p\,= \,\bigl\{\,x\in
l_{p}\,\bigl| \, x \ge 0,\, \sum_{i=1}^{+\infty}\, x^i\, \le\, 1\,
\bigr\}$ in the space $l_{p}$ is pointwise $\mu$-compact, but not
$\mu$-compact.}
\end{property}
\smallskip

\noindent {\tt Proof.} The pointwise $\mu$-compactness of the set
$\Delta_{p}$ follows from the remark above. Let us show that
$\Delta_{p}$ is not $\mu$-compact in the space~$l_{p}$. Consider an
increasing sequence of natural numbers $\{n_r\}_{r \in \mathbb N}$
and a sequence of nonnegative numbers $\{z_i\}_{i \in \mathbb N}$
such that $\, \sum_{i=n_{r}}^{n_{r+1}-1}\, z_i \, = \, 1\, $ and $\,
\sum_{i=n_{r}}^{n_{r+1}-1}\, (z_i)^p \, \le\, \frac1r\, $ for each
$\, r \ge 1$. The set $\, \mathcal{K}\,= \, \bigl\{ \, y \in
\Delta_p \, \bigl| \, \forall \, r \in \mathbb N \
\sum_{i=n_r}^{+\infty}\, (y^i)^p \, \le  \, \frac1r \bigr\}$ is
compact. If the set~$\Delta_p$ is $\mu$-compact, then by
Proposition~\ref{ce-prop-0} for any $\varepsilon > 0$
 there is the corresponding compact set~$\, \mathcal{K}_{\varepsilon}$.
The set $\mathcal{K}_{\varepsilon}$ can contain only a finite number
of vectors of the canonical basis $\{e_i\}$, otherwise it is not
compact. Let $N$ be such that $\, e_i \, \notin \,
\mathcal{K}_{\varepsilon}$ for $i
> N$. Take an arbitrary $r$ for which $n_r > N$ and denote by
$x \in \Delta_p$ the vector with coordinates $\, x^i \, = \, z^i\, $
for $\, n_r \, \le \, i \, <
 \, n_{r+1}\, $ and $\, x^i = 0$ for all other $i$. Clearly, $x \in \mathcal{K}$.
Since
 $x = \sum_{i=n_r}^{n_{r+1}-1}\, x^i\, e_i\, $ and $\, \sum_{i=n_r}^{n_{r+1}-1} x^i = 1$,
but  $\, e_i \notin \mathcal{K}_{\varepsilon}\, $ for all $i \,=\,
\overline{n_r,n_{r+1}-1}$, we come to the contradiction with the
$\mu$-compactness criterion from Proposition~\ref{ce-prop-0}.

{\hfill $\Box$}
\medskip

In particular, Proposition~\ref{l1-lp} shows that a Hilbert space
contains noncompact pointwise  $\, \mu$-compact sets. It appears
that it does not contain $\, \mu$-compact sets that are not compact.

\medskip

\begin{property}\label{l2} {\it There are no  $\, \mu$-compact subsets of Hilbert space that are not compact.}
\end{property}
\smallskip

\noindent {\tt Proof}. We describe the idea of the proof omitting
details, which can be easily reconstructed by the reader. Let
$\mathcal{A}$ be a bounded convex closed subset of a Hilbert space.
Without loss of generality it can be assumed that its diameter
is~$1$. If~$\mathcal{A}\, $ is not compact, then there exists a
sequence of its elements  $\{\, a_k\, \}_{\, k \, \in \, \mathbb
N}\, $ such that all their norms and all pairwise distances between
them exceed some value~$\, \varepsilon \, > \, 0$. Because of the
weak compactness of~$\mathcal{A}$ it can be assumed, after possible
passing to a subsequence, that the sequence  $\,\{a_k\}\, $
converges weakly to some element $\, a$. Since the set $\mathcal{A}$
is convex and closed, we have $\, a \, \in \, \mathcal{A}\, $. The
sequence $\,\{ b_{\, k}\, = \, a_{\, k}\, - \, a\,\}_{\, k \, \in \,
\mathbb N}\,$ converges weakly to zero. Hence by passing to
subsequences one can achieve a fast convergence to zero for the
scalar products:
  $\, (b_{\, i}, b_{\, j})\, \le \, \varepsilon^{\, 2}\, 2^{\, -\, i\, -\, j}\, $
 for every $\ i\, \ne \, j$.
Since ${\, \varepsilon \, \le \, \|\, b_{\, i}\, \|_{\, 2}\, \le \,
1\, }$ for all $\,  i \, \in \, \mathbb N\, $,  invoking the
Cauchy-Schwarz inequality we conclude that any sequence $\, x \, \,
\in \, l_{\, 2}\, $ satisfies the inequalities
$$
{\ \frac{\varepsilon}{2}\, \bigl\|\,  x  \, \bigr\|_{\, 2} \ \le \
\Bigl\|\, \sum_{i}\, x^{\, i}\, b_{\, i}    \, \Bigr\|_{\, 2}    \
\le \     2\,  \bigl\|\,  x \, \bigr\|_{\, 2} }.
$$
This means that the system of elements  $\, \{b_k\}\, $ possesses
the Riesz basis property. Then there is a continuous linear operator
that is continuously invertible and taking the system
 $\, \{\, b_k \, \}\, $ to an orthonormal one~\cite{KS}. This operator maps the convex hull
 of points  $\, a\, $ and $\, \{\, a_{k}\, \}_{\, k \,
\in \, \mathbb N}\, $ to the set $\, \Delta_{\, 2}\, $ from
Proposition~\ref{l1-lp}, which is not $\, \mu$-compact. Hence, this
convex hull is not $\, \mu$-compact, therefore nor is the set~$\,
\mathcal{A}\, $.

{\hfill $\Box$}
\medskip

Let us now consider two important examples of not pointwise
$\mu$-compact sets.
\medskip

\begin{property}\label{lp}
{\em For every $\, p > 1$ the set  $\, \mathcal{A}_p \, = \,
\bigl\{\, x \in l_p\;|\, \ x \ge 0,\ \|x\|_p \, \le \, 1\, \bigr\}$
(the bounded part of the positive cone in the space~$l_p$) is not
pointwise $\mu$-compact}.
\end{property}
\smallskip

\noindent {\tt Proof.} Let us show that if a point $\, x \, \in \,
\mathcal{A}_p\, $ is such that $\, \, \| x\|_p\, < \, \frac13$ and
$\sum_i x_i\, = +\infty$, then it has no compact set $\,
\mathcal{K}_{\varepsilon}$ for ${\, \varepsilon \,=\, \frac13\, }$
(see Proposition~\ref{ce-prop-0}). If such a compact set exists,
then it can contain only finitely many elements of the canonical
basis~$\{e_i\}$. Take a sufficiently  large  $N$ such that
 $\, e_i \, \notin \, \mathcal{K}_{\varepsilon}$ for every $i > N$.
Since the series $\, \sum_i x_i\, $ diverges, we see that there
exists $r$, for which $\, s \,  = \, \sum_{i=N+1}^{N+r} x^i\, \in \,
\bigl(\frac13 \, , \, \frac23\bigr)$. Let $\, \bar x \, = \, x \, -
\, \sum_{i=N+1}^{N+r} x^i\, e_i$ (the entries of the vector $\bar
x\, $ from the $\, (N+1)$st to the $(N+r)$th  are zeros, and the
others coincide with the corresponding entries of~$x$). Then $\ x
\,= \, (1-s)\, \bigl( \frac{1}{1-s}\, \bar x\bigr) \, + \,
\sum_{i=N+1}^{N+r} x^ie_i$. Since $\, \frac{1}{1-s} \, < \, 3$ and
$\, \|\bar x\|_p \, < \, \frac13$, it follows that $\,
\frac{1}{1-s}\, \bar x \, \in \, \mathcal{A}_p$. All the
points~$e_i$ in this barycentric  combination lie outside of
$\mathcal{K}_{\varepsilon}$, but their total weight~$s$ exceeds
$\frac{1}{3}$, which is a contradiction.$\;\Box$
\medskip

The next example complements Proposition~\ref{formula-} and
demonstrates that not pointwise
$\mu$\nobreakdash-\hspace{0pt}compact sets do not have to satisfy
relations (\ref{co-cl-exp}) and (\ref{co-cl-exp+}) for continuous
bounded concave functions. This violates
representation~(\ref{ext-rep+}).\medskip

\begin{example}\label{conv-cl-ex}
Let $f$ be a continuous function on the bounded part $\mathcal{A}_p$
of the positive cone of the space~$l_{p}$ for $p>1$ that takes the
value~$1$ at zero and vanishes at all the vectors of the canonical
basis~$\{e_{n}\}$ of the space $l_{p}$. For example, the function
$f(\cdot)=1-\|\cdot\|_{p}$. Since the zero vector of the
space~$l_{p}$ is a limit point of the set of all convex combinations
of the vectors~$\{e_{n}\}$, the convex closure of the function~$f$
is an identical zero on the set~$\mathcal{A}_p$. Therefore,
$\overline{\mathrm{co}}f(0)\neq f(0)$. This example shows also that
for the set~$\mathcal{A}_p$ such that
$\mathcal{A}_p=\overline{\mathrm{co}}(\mathrm{extr}\mathcal{A}_p)$
and map~(\ref{b-map}) is open (this follows from a theorem
in~\cite{Grzaslewicz} combined with the strict convexity of the
space~$l_{p}$ for $p>1$ and with the proof of theorem~1
in~\cite{Sh-8}) the assertion of corollary~2 in~\cite{Sh-8} does not
hold, that is, a continuous function on the closed
set~$\mathrm{extr}\,\mathcal{A}_p$ need not have a convex continuous
(or even lower semicontinuous) extension to~$\mathcal{A}_p$.
\end{example}

A Hilbert cube is the set $\, \mathcal{H}_a\, =\, \bigl\{\, x \, \in
\, l_2 \ \bigl| \ |\, x^{\, i}\, | \, \le \, a^{\, i}\, , \ i \in
\mathbb N\, \bigr\}$, where $\, a \, =\, \{\, a^{\, i}\}_{\, i \in
\mathbb N}$ is an arbitrary sequence of positive numbers. Let us
show that the following alternative holds for any Hilbert cube.
\smallskip

\begin{property}\label{brick}
\textit{If $\; \|\, a \, \|_2 <  +\infty$, then the set $\,
\mathcal{H}_a \, $ is compact; otherwise, if ${\; \|\, a \,
\|_2=+\infty }$, then it is not even pointwise $\mu$-compact.}
\end{property}
\smallskip

\noindent {\tt Proof}. The first assertion is well known. It
follows, for instance, from the compactness criterion for~$l_2$. If
${\, \|\, a \, \|_2\, = +\infty }$, then we split the sequence $\,
a\, $ into blocks so that each block contains finitely many
consecutive elements the sum of whose squares exceeds~$1$.  Suppose
that the $\, n$\nobreakdash-\hspace{0pt}th block consists of
elements $\, a^{\, k}, \,  a^{\, k+1}\, , \,  \ldots \ , \, a^{\,
k+m}$. We write $\, b_{\, n} \, = \, \sum_{i \, = \, k}^{\, k \, +
\, m}\, a^{\, i}\, e_i$ for every $n$ and denote by $\mathcal{L}$
the closure of the linear span of the elements $\, b_{\, n}\, , \, n
\in \mathbb N$. The set~$\mathcal{L}$ is a Hilbert space with
orthogonal basis~$\, \{b_{\, n}\}$. The unit ball of the
space~$\mathcal{L}$ is contained in $\, \mathcal{H}_a$. Since a ball
in a Hilbert space is not pointwise $\mu$-compact, it follows that
the set $\, \mathcal{H}_a\, $ is not so either.

{\hfill $\Box$}
\medskip

\section{The $\mathrm{CE}$-properties of $\mu$-compact convex sets}

In the 1970s various  properties of convex compact sets, in
particular, the continuity properties of the convex hulls of
continuous functions have been  intensively studied. Vesterstrom
proved in \cite{Ves} the relation between continuity of the convex
hull \footnote{By corollary I.3.6 in \cite{Alf} the convex hull of
any continuous function on a compact set coincides with the convex
closure of this function.} of arbitrary continuous function and the
openness of the barycenter map. He also conjectured the equivalence
between continuity of the convex hull of any continuous concave
function (called the $\mathrm{CE}$\nobreakdash-\hspace{0pt}property
in \cite{Lima}) and the continuity of the convex hull of any
continuous function (called the strong
$\mathrm{CE}$\nobreakdash-\hspace{0pt}property in \cite{Sh-8}). This
conjecture was approved by O'Brien \cite{Brien}, who showed the
equivalence of these properties to the openness of the convex
mixture map $(x, y, \lambda)\mapsto\lambda x+(1-\lambda)y$. In the
subsequent papers the latter property was studied for convex sets
that are not necessarily compact, and was called the
\textit{stability} property. Convex sets having this property were
called \textit{stable} convex sets \cite{Susanne}. The relations
between the stability property and several other properties of
convex sets were also revealed \cite{Clausing},\cite{Grzaslewicz}.

In this section we  generalize the Vesterstrom-O'Brien theory to the
class of $\mu$\nobreakdash-\hspace{0pt}compact convex sets. The
first partial result in this direction was obtained in \cite{Sh-8},
where the $\mu$\nobreakdash-\hspace{0pt}compact version of theorem
3.1 in \cite{Ves} was proved. The following theorem is the
$\mu$\nobreakdash-\hspace{0pt}compact generalization of the main
result in \cite{Brien}.\vspace{5pt}

\begin{theorem}\label{main}
\textit{For a convex $\mu$-compact set $\mathcal{A}$ the following
properties are equivalent:}
\begin{enumerate}[(i)]
 \item \textit{the map $\ \mathcal{A}\times\mathcal{A}\, \ni\, (x,y)\ \mapsto
  \, \displaystyle\frac{x\, +\, y}{2}\, \in\, \mathcal{A}$ is open (the stability property);}
 \item \textit{the map $M(\mathcal{A})\ni\mu\ \mapsto \ \mathbf b(\mu)\in \mathcal{A}$
is open;}
 \item \textit{the map $M\left(\overline{\mathrm{extr}\,
\mathcal{A}}\right)\ni\mu\ \mapsto \
 \mathbf b(\mu)\in \mathcal{A}$
is open;}\footnote{This map is surjective by Proposition
\ref{KMC}.}
 \item \textit{the convex hull of an arbitrary function in $C(\mathcal{A})$ is continuous  (the strong CE-property);}
 \item \textit{the  convex hull of an arbitrary function in $Q(\mathcal{A})$ is continuous  (the CE-property).}
\end{enumerate}
\smallskip

\textit{Equivalent properties  $\mathrm{(i)-(v)}$ imply the
closedness of the set~$\mathrm{extr}\, \mathcal{A}$.}
\end{theorem}
\smallskip

\begin{remark}\label{main+}
Property $\mathrm{(i)}$ in Theorem~\ref{main} is equivalent to the
openness of the map \break
$\mathcal{A}\times\mathcal{A}\times[0,1]\, \ni\, (x,y,\lambda)\
\mapsto\ \lambda \, x\, +\, (1-\lambda)\, y\, \in\, \mathcal{A}\;$
\cite{Clausing}. Properties $\mathrm{(iv)}$ and $\mathrm{(v)}$ in
Theorem~\ref{main} can be formulated as the continuity of the
convex closure and its coincidence with the convex hull for any
function in $C(\mathcal{A})$ and in  $Q(\mathcal{A})$
respectively.
\end{remark}

\begin{remark}\label{main++}
If properties $\mathrm{(i)-(v)}$ hold for $\mu$-compact set
$\mathcal{A}$, then the family $F(\mathcal{A})$ in
Proposition~\ref{ce-prop-1} can be chosen consisting of lower
semicontinuous functions. Indeed, by using property
$\mathrm{(ii)}$ it is easy to show that the functions
$f_{\varphi}$ constructed in the proof of
Proposition~\ref{ce-prop-1} are lower semicontinuous.
\end{remark}\vspace{5pt}

{\tt The proof of Theorem~\ref{main}.} Note first that
$\mathrm{(v)}$ implies the closedness of the set $\mathrm{extr}\,
\mathcal{A}$  by Proposition~\ref{ext-rep}, since $\mathrm{(v)}$
guarantees the closedness of the set $\mathcal{B}_{f}=\{x\in
\mathcal{A}\,|\,f(x)=\overline{\mathrm{co}}f(x)\}$ for any function
$f\in Q(\mathcal{A})$.

$\mathrm{(v)}\Rightarrow\mathrm{(iii)}$ The proof of this part of
theorem (as well as the proof of the analogous part of theorem 3.2
in \cite{Ves}) can be realized by means of lemma 2.1 in
\cite{Ves}. That lemma can be proved  without the compactness
assumption by  the following observation: \textit{if $X$ is a
compact space and $Y$ is an arbitrary topological space then the
image of any closed subset of $X\times Y$ under the canonical
projection $X\times Y\ni(x,y)\mapsto y\in Y$ is a closed subset of
$Y$}.

By representation (\ref{co-cl-exp}) the convex closure of any
function $f$ in $Q(\mathcal{A})$ is determined by the expression
$$
\overline{\mathrm{co}}\, f(x)\ =\ \inf_{\mu\in
M_{x}(\mathrm{extr}\mathcal{A})}\mu(f),\quad \forall x\in
\mathcal{A},\quad \textrm{where} \quad \mu(f)\ =\
\int_{\mathrm{extr}\mathcal{A}} f(x)\, \mu(dx).
$$
Hence $\mathrm{(v)}$ yields  the continuity and boundedness of the
function
$$
\mathcal{A}\ni x\mapsto \sup\{\mu(f)\,|\,\mu\in M(\mathrm{extr}\,
\mathcal{A}), b(\mu)=x\}
$$
for any function $f$ in $P(\mathcal{A})$. The above-mentioned
generalization of lemma 2.1 in \cite{Ves} with $K=M(\mathcal{A})$,
$M=M(\mathrm{extr}\mathcal{A})$ and $K'=\mathcal{A}$ implies the
openness of the map
\begin{equation}\label{b-map+}
M(\mathrm{extr}\,\mathcal{A})\ \ni \ \mu\quad \mapsto \quad \mathbf
b(\mu)\ \in \ \mathcal{A},
\end{equation}
provided the set $M(\mathrm{extr}\, \mathcal{A})$ is endowed with
the topology, which has the prebase consisting of the sets $\{\mu\in
M(\mathrm{extr}\, \mathcal{A})\,|\,\mu(f)>0\},\; f\in
P(\mathcal{A})$. Following the terminology in \cite{Ves}, this
topology will be called the $p$-topology. This is the weakest
topology providing lower semicontinuity of the functionals
$\mu\mapsto\mu(f)$ for any function $f\in P(\mathcal{A})$.

By using Lemma~\ref{open-map} (see Appendix) we will show that the
openness of map (\ref{b-map+}) in the
$p$\nobreakdash-\hspace{0pt}topology on $M(\mathrm{extr}\,
\mathcal{A})$ and the closedness of the set $\mathrm{extr}\,
\mathcal{A}$ proved above imply the openness of map (\ref{b-map+})
in the weak topology on $M(\mathrm{extr}\,
\mathcal{A})$.\footnote{In the case of compact set $\mathcal{A}$ the
coincidence of these topologies on $M(\mathrm{extr}\, \mathcal{A})$
is proved in \cite{Ves} (lemma 3.4). In the case of $\mu$-compact
set $\mathcal{A}$ we can not prove this coincidence.} To realize
this it suffices to show that for an arbitrary converging sequence
$\{x_{n}\}\subset\mathcal{A}$ and for an  arbitrary net
$\{\mu_{\lambda}\}_{\lambda\in\Lambda}\subset M(\mathrm{extr}\,
\mathcal{A})$ such that
$$
b(\{\mu_{\lambda}\}_{\lambda\in\Lambda})\subseteq\{x_{n}\}\quad\textrm{and}\quad\exists\;
p\,\textrm{-}\lim_{\lambda}\mu_{\lambda}\ =\ \mu_{0},
$$
where $\mu_{0}$ is a measure in $M(\mathrm{extr}\, \mathcal{A})$
such that $\mathbf b(\mu_{0})=\lim_{n\rightarrow+\infty}x_{n}$,
there exists a subnet of the net
$\{\mu_{\lambda}\}_{\lambda\in\Lambda}$ that   converges weakly to
the measure $\mu_{0}$.

Let $\{x_{n}\}$ and $\{\mu_{\lambda}\}_{\lambda\in\Lambda}$ be the
above sequence and net, respectively. Since the sequence is
relatively compact, $\mu$-compactness of the set $\mathcal{A}$ and
the inclusion $\mathbf
b(\{\mu_{\lambda}\}_{\lambda\in\Lambda})\subseteq\{x_{n}\}$ imply
relative compactness of the net
$\{\mu_{\lambda}\}_{\lambda\in\Lambda}$ in the weak topology and
hence the existence of the subnet
$\{\mu_{\lambda_{\pi}}\}_{\pi\in\Pi}$ weakly converging to a some
measure $\nu\in M(\mathrm{extr}\mathcal{A})$. By the definitions of
the weak topology and of the $p$-topology we have
$$
\nu(f)\ =\ \lim_{\pi}\mu_{\lambda_{\pi}}(f)\ \geq\
p\,\textrm{-}\liminf_{\lambda}\, \mu_{\lambda}(f)\ \geq \
\mu_{0}(f),\quad \forall \, f\, \in \, P(\mathcal{A}).
$$
This means that $\nu\succ\mu_{0}$ (in the Choquet ordering).
Closedness of the set $\mathrm{extr}\mathcal{A}$ implies the
maximality in $M(\mathcal{A})$ of any measure in
$M(\mathrm{extr}\mathcal{A})$. This can be proved by using Theorem
2.2 in \cite{Edgar-2} and the arguments from the proof of Theorem
1.1 in \cite{B&E}, but it can also be immediately shown by using
property $(\mathrm{v})$ and coincidence of any function in
$Q(\mathcal{A})$ with its convex hull on the set
$\mathrm{extr}\mathcal{A}$. Thus $\mu_{0}$ is the maximal measure in
$M(\mathcal{A})$ and hence $\nu=\mu_{0}$.

$\mathrm{(iii)}\Rightarrow\mathrm{(i)}\ $ This implication follows
from Proposition~\ref{KMC} and Proposition~\ref{last-prop} below
(with $X=\overline{\mathrm{extr}\mathcal{A}}$).
\smallskip

By Remark~\ref{main+} the equivalence of properties
$\mathrm{(i)}$, $\mathrm{(ii)}$ and $\mathrm{(iv)}$ for convex
$\mu$-compact subsets of a Banach space is proved in \cite[theorem
1]{Sh-8}. That proof is easily extended  to the class of sets
considered in this paper.

The implication $\mathrm{(iv)}\Rightarrow\mathrm{(v)}$ is obvious.

{\hfill $\Box$}
\medskip

In the proof of Theorem~\ref{main} we have involved the following
result of the measure theory (\cite[theorem 2.4]{Eifler}).
\smallskip

\begin{property}\label{last-prop}
\textit{Let $X$ be a complete separable metric space. The map
$$
M(X)\times M(X)\ni(\mu,\nu)\mapsto\frac{1}{2}(\mu+\nu)\in M(X)
$$
is open.}
\end{property}
\smallskip

\noindent Since the set $M(X)$ is $\,\mu$-compact
(Corollary~\ref{p-measures}), we arrive at the following
observation.
\medskip

\begin{corollary}\label{p-measures+}
\textit{Properties  $\mathrm{(i)} - \mathrm{(v)}\ $ in
Theorem~\ref{main} hold for the set of Borel probability measures
on a complete separable metric space endowed with the weak
convergence topology.}
\end{corollary}
\medskip

The $\mu$-compactness condition is essential in the proof of
Theorem~\ref{main}, it can not be removed without changing the
whole structure of the proof. This motivates the conjecture that
the class of $\mu$-compact convex sets is the maximal class of
convex metrizable sets for which the Vesterstrom - O'Brien theory
can be generalized. This conjecture can be justified by the
following example, showing that even pointwise $\mu$-compactness
is not sufficient for the proof of Theorem~\ref{main}.
\medskip

\begin{property}\label{counter-example}
\textit{For any $p>1$ the pointwise $\mu$-compact simplex
$$
\Delta_p\, = \, \bigl\{\,x\in l_{p}\ \bigl| \  x\,  \ge \, 0\, ,\
\sum_{i=1}^{+\infty}x^i\, \le \, 1\, \bigr\}
$$
in $l_{p}$ is stable, that is, it possesses property $\mathrm{(i)}$
in Theorem~\ref{main}, which implies $\mathrm{(ii)}$, but it does
not possess properties $\mathrm{(iii)}-\mathrm{(v)}$.}
\end{property}
\medskip

Note that for $p=1$ the $\mu$-compact simplex
~$\Delta_1=\mathcal{A}_1$ has properties $\mathrm{(i)} -
\mathrm{(v)}$ in Theorem~\ref{main}.
\medskip

\noindent {\tt Proof.} Example \ref{conv-cl-ex} shows that the
simplex $\Delta_p$ does not possess properties $\mathrm{(iv)}$ and
$\mathrm{(v)}$. Let us show that it does not possess
$\mathrm{(iii)}$ either. Note that $\, {\rm extr}\, \Delta_p\, = \,
\bigl\{\, 0\, , \, e_i \, , \ i \, \in \, \mathbb N \, \bigr\}$ and
that the set $\, \Delta_p\, $ is a real simplex: for any its point
$x$ there exists the unique measure on $\, {\rm extr}\, \Delta_p\, $
with the barycenter~$x$. The sequence of points $\, x_n \, = \,
\bigl( \frac1n \, , \, \cdots \, \frac1n , 0 \, , \cdots \bigr) \,
\in \, \Delta_p$ (the first $n$ coordinates equal to $\, \frac1n$,
the all others are zeros) converges to zero in $l_p$ as $\, n \, \to
\, +\infty$, but it is easy to see that the corresponding sequence
of measures on $\, {\rm extr}\, \Delta_p\, $ does not converge to
the single atom measure supported at the point $0$.

Let us now show that for an arbitrary $\, p \, > \, 1\, $ the set
$\, \Delta_p\, $ is stable, that is,  it has property
$\mathrm{(i)}$, which implies  $\mathrm{(ii)}$ (see the proof of
theorem 1 in \cite{Sh-8}). It suffices to prove that for arbitrary
points $\, a , b \, \subset \, \Delta_p \, , \ c\, = \,
\frac{1}{2}\bigl(\, a\, + \, b\, \bigr)\, $ and for arbitrary $\,
\varepsilon \, > \, 0\, , $ there exists $\,\delta \, > \, 0\, $
with the following property: for any $\, z \, \in \, \Delta_p\, $
such that $\, \|\, z \, - \, c \, \|_p\, < \, \delta \, $  there
exists a segment $\, [x,y]\, \subset \, \Delta_p \, $ centered at
$z$, for which $\, \|\, x \, - \, a \, \|_p\, < \, \varepsilon \, $
and  $\, \|\, y \, - \, b \, \|_p\, < \, \varepsilon $. By taking
sufficiently small $\, \varepsilon\, $ it can be assumed that $\,
\|\, a\, \|_p \, < \, 1\, -\, \varepsilon \, $ and that $\, \|\, b\,
\|_p \, < \, 1\, -\, \varepsilon$. Otherwise the points $\, a\, $
and $\, b\, $ can be replaced by sufficiently close points belonging
to the interior of $\, [a\, ,\, b]$. Since the $l_p$-norm is
strictly convex, the norms of $\, a\, $ and of $\, b\, $ become less
then $1$. Then we choose a large~$N$ so that for each of the points
$\, a\, $ and $\, b\, $ the norm of the ''tail'' from the $(N+1)$th
coordinate is less than $\frac16 \, \varepsilon$, that is, $\,
\bigl(\sum_{k=N+1}^{+\infty} \, (a^k)^p\bigr)^{1/p}\,  < \, \frac16
\, \varepsilon\, $ and $\, \bigl(\sum_{k=N+1}^{+\infty} \,
(b^k)^p\bigr)^{1/p}\,  < \, \frac16 \, \varepsilon\, $. Consider the
space $\, \mathbb R^N$, generated by the first $N$ coordinates.
Denote by $\, \tilde \Delta_p\, $ and  $\, \tilde s$ the
restrictions of the set $\, \Delta_p\, $ and of arbitrary element
$\, s \, \in \, l_p\, $ to this space. Since $\, \tilde \Delta_p\, $
is a simplex in $\, \mathbb R^N\, $ it is stable (see
\cite{Susanne}) and hence one can take  $\, \delta \, >\, 0\, $ such
that there exist points $\, \tilde x\, , \, \tilde y \, \in \,
\tilde \Delta_p,  \, $ for which $\, \frac{1}{2}\, \bigl( \, \tilde
x \, + \, \tilde y \, \bigr)\,
 = \, \tilde z \, , $   $\, \|\, \tilde x \, - \, \tilde a  \,  \|_p \, < \,
\frac13 \, \varepsilon\, $ and  $\, \|\, \tilde y \, - \, \tilde
b\, \|_p\, < \, \frac13\, \varepsilon  \, , $ if  $\, \|\, \tilde
z \, - \,  \tilde c\, \|_p \,  < \, \delta$. Now for given $\, z
\, \in \, \Delta_p \, $ and for arbitrary  $t \in [-1, 1]$ we
define the points $\, x(t)  ,  \, y(t) \, \in \, l_p\, $ as
follows:
$$
x^{\, k}(t)\ =\ \left\{
\begin{array}{lcl}
\ \tilde x^{\, k}\,  , & k \, \le \, N \\ (1+t)\, z^{\, k} \, , & k
\, > \, N
\end{array}
\right.   \qquad y^{\, k}(t)\ =\ \left\{
\begin{array}{ll}
\ \tilde y^{\, k}\,  , & k \, \le \, N \\ (1-t)\, z^{\, k} \, , & k
\, > \, N
\end{array}
\right.
$$
By construction $\, \frac12 \, \bigl( \, x(t)  \, + \, \, y(t)\,
\bigr)\, = \, z\, $ for any $t$, while the norms of the elements
$x(t)\, $ and  $\, y(t)\, $ do not exceed  $\,  1 \, -\, \frac13\,
 \varepsilon \, + \,
2\, \delta$. Indeed,
$$
 \bigl\|\, \tilde x\, \bigr\|_p\ \le \ \bigl\|\, \tilde a\, \bigr\|_p   \ + \  \bigl\|\,
 \tilde x\, - \tilde a\, \bigr\|_p
\ \le \ 1 \ - \ \varepsilon \ + \ \frac13 \, \varepsilon \ = \ 1 \ -
\  \frac23 \, \varepsilon,
$$
while the norm of the ``tail'' of $\, (1+t)\, z \, $ does not
exceed the sum of norms of the ``tails'' of
 $\, 2 c \, $ and ${\, 2\, (z-c)}$, that is, it does not exceed
 $\, 2\, \bigl( \,  \frac16\, \varepsilon \, + \, \delta   \, \bigr)$.
Taking the sum we obtain  $\, \bigl\|\, x(t)\, \bigr\|_p \, \le \, 1
\, -\,
 \frac13\, \varepsilon \, + \,  2\, \delta $, and the same is
 true for $\, y(t)$. For  $\, \delta \, \le \, \frac{1}{6}\, \varepsilon \, $
we get $\, \bigl\|\, x(t)\, \bigr\|_p \, \le \, 1 \ $ and
$\bigl\|\, y(t)\, \bigr\|_p \, \le \, 1 \, $. Let us show that
there exists $\, \tau \, \in \, [-1,1]$, for which $\, \| \, x\,
(\tau)\, \|_1 \, \le \, 1\, $ and  $\, \|\,  y(\tau)\, \|_1 \, \le
\, 1$, and whence ${\, x\, (\tau),\, y\, (\tau)\, \in \, \Delta_p
\, }$. It is clear that $\, \|\, \tilde x\, \|_1 \, \le \, 1\, $
and $\, \|\, \tilde y\, \|_1 \, \le \, 1$. For the sake of being
defined, suppose  $\, \|\, \tilde x\, \|_1 \, \ge \, \|\, \tilde
y\, \|_1$. If ${\, \|\, y(-1)\, \|_1\, \le \, 1}$, then one can
take $\, t \, = \, -1$, since $\, \|\, x(-1)\, \|_1\, = \, \|\,
\tilde x\, \|_1\, \le \, 1$. If $\, \|\, y(-1)\, \|_1\,
> \, 1$, then $\, \|\, y(-1)\, \|_1\, > \, \|\, x(-1)\, \|_1$, and since
 $\, \|\, y(1)\, \|_1\,  \le  \, \|\, x(1)\|_1$, applying the continuity argument, we conclude
 that there is  $\, \tau \, \in \,
[-1,1]$ such that   $\, \|\, y(\tau)\, \|_1\, = \, \|\, x(\tau)\,
\|_1$. Since $\, \frac12 \, \bigl( \, x(\tau)  \, + \,  \,
y(\tau)\, \bigr)\, = \, z\, $, we have $\, \|\, x(\tau)\, \|_1\, =
\, \|\, y(\tau)\, \|_1\, = \, \|\, z\, \|_1 \, \le \, 1$.

Finally, the norm of the difference  $\, \|\, x(\tau) \, - \, a\,
\|_p\, $  does not exceed $\, \frac13\, \varepsilon$ for the first
$N$ coordinates, while for the other coordinates it does not exceed
the maximal norm of the two ``tails'': of the element $\, a\, $ and
of~$\, 2 z$. Hence,
$$
 \bigl\|\, x(\tau) \ - \ a\, \|_p\ \le \ \frac13\, \varepsilon \ + \
\max \Bigl\{ \,  \frac16 \, \varepsilon \ , \ 2\, \Bigl(\, \frac16
\, \varepsilon  \ + \ \delta \, \Bigr)     \, \Bigr\}\ = \ \frac23
\, \varepsilon \  + 2\, \delta \, .
$$
For $\, \delta  \, < \, \frac16\, \varepsilon $ we obtain $\, \|\,
x(\tau) \, - \, a\, \|_p \, <  \, \varepsilon $, and similarly $\,
\|\, y(\tau) \, - \, b\, \|_p \, < \, \varepsilon $. By setting
$\, x \, = \, x(\tau)\,  ,  \, y \, = \, y(\tau)$, we complete the
proof.

{\hfill $\Box$}

\medskip

\section{Applications to quantum information theory}

An important example of convex $\mu$-compact sets, for which the
equivalent properties of Theorem~\ref{main} hold, is the set
$\mathfrak{S}(\mathcal{H})$ of quantum states.\footnote{The set
$\mathfrak{S}(\mathcal{H})$ is compact if and only if
$\dim\mathcal{H}<+\infty$.} The quantum states are  density
operators (positive operators with trace equal to $1$) in a
separable Hilbert space $\mathcal{H}$ \cite{H-QI}. The set of
extreme points of the set $\mathfrak{S}(\mathcal{H})$ consists of
one-dimensional projectors called pure states. The $\mu$-compactness
and the stability property $\mathrm{(i)}$ in Theorem~\ref{main} of
the set $\mathfrak{S}(\mathcal{H})$ are established  in
\cite[proposition 2]{H-Sh-2} and in \cite[lemma 3]{Sh-2}
respectively. These properties are  used essentially in the study of
characteristics of quantum states and of quantum channels. For
instance, $\mu$-compactness of the set $\mathfrak{S}(\mathcal{H})$
makes it possible to prove that any nonentangled state (see below)
of composite quantum system can be represented as an average
(barycenter)  state
 of some generalized ensemble of pure product states
(probability measure on the set of product pure states)
\cite{H-Sh-W}. The stability property of the set
$\mathfrak{S}(\mathcal{H})$ plays the crucial role  in the proof of
the lower semicontinuity property of the $\chi$-function of an
arbitrary quantum channel, which is an important characteristic
related to the classical capacity of this channel \cite{Sh-2}.

In this section we consider a result  following directly from
 the generalized Vesterstrom-O'Brien theorem (Theorem~1)
applied to the stable set $\mathfrak{S}(\mathcal{H})$.

According to the quantum mechanical formalism, states of composite
quantum system, arising as a result of joining  two quantum systems
represented by two  Hilbert spaces $\mathcal{H}$ and $\mathcal{K}$,
correspond to the density operators in the tensor product
$\mathcal{H}\otimes\mathcal{K}$ of these spaces. A specific property
of the quantum mechanical statistical model (in comparison to the
classical one) is the existence of the so called \textit{entangled}
states of composite system, which can not be represented as convex
combinations of product states describing independent subsystems.
Entanglement can be considered as a special purely quantum
correlation, which is the base
 for the construction of different quantum
algorithms, quantum criptographical protocols and systems of
information transmissions, which attract a lot of attention of
scientists in the last two decades (see \cite{H-QI}, chapter 3).
That is why the study of the entanglement, in particular, of  its
quantitative characteristics, is one of the main problems of the
quantum information theory.

Let $\mathcal{H}$ and $\mathcal{K}$ be separable Hilbert spaces. A
state $\omega\in\mathfrak{S}(\mathcal{H}\otimes\mathcal{K})$ is
called \emph{nonentangled} if it belongs to the convex closure of
the set of product states, that is, of states of the form
$\rho\otimes\sigma$, where $\rho\in\mathfrak{S}(\mathcal{H})$ and
$\sigma\in\mathfrak{S}(\mathcal{K})$; otherwise it is called
\emph{entangled}.

\textit{Entanglement monotone} is an arbitrary function $E$ on the
set $\mathfrak{S}(\mathcal{H}\otimes\mathcal{K})$ that possesses
the following properties (see \cite{Vidal},\cite{P&V}):

E-1)
$\{E(\omega)=0\}\Leftrightarrow\{\textit{the}\;\textit{state}\;\omega\;\textit{is}\;\textit{nonentangled}\}$;

E-2) \textit{Monotonicity under Local Operations and Classical
Communications (LOCC),} which means
\begin{equation}\label{LOCC-m}
E(\omega)\geq\sum_{i}\pi_{i}E(\omega_{i})
\end{equation}
for arbitrary state
$\omega\in\mathfrak{S}(\mathcal{H}\otimes\mathcal{K})$ and arbitrary
LOCC-operation, transforming the state $\omega$ into the set
$\{\omega_{i}\}$ of states with the probability distribution
$\{\pi_{i}\}$ (see details in \cite{P&V}).

E-3) \textit{Convexity of the function $E$ on the set
$\mathfrak{S}(\mathcal{H}\otimes\mathcal{K})$}, which means
$$
E\left(\sum_{i}\pi_{i}\omega_{i}\right)\leq
\sum_{i}\pi_{i}E(\omega_{i})
$$
for arbitrary finite set $\{\omega_{i}\}$ of states in
$\mathfrak{S}(\mathcal{H}\otimes\mathcal{K})$ and probability
distribution  $\{\pi_{i}\}$.

The standard method of "generation" of entanglement monotones (EM)
in the case of finite dimensional spaces $\mathcal{H}$ and
$\mathcal{K}$ is the convex roof construction (see
\cite{P&V},\cite{Osborne}). In accordance with this method, for an
arbitrary concave continuous nonnegative function $f$ on the set
$\mathfrak{S}(\mathcal{H})$ such that
\begin{equation}\label{f-property}
f^{-1}(0)=\mathrm{extr}\mathfrak{S}(\mathcal{H})\quad\textup{and}\quad
f(\rho)=f(U\rho U^{*})
\end{equation}
for any state $\rho $ in $\mathfrak{S}(\mathcal{H})$ and any unitary
operator $U$ in the space $\mathcal{H}$ the corresponding EM $E^{f}$
is defined by
\begin{equation}\label{E-f}
E^{f}(\omega)=\inf_{\{\pi_{i},\omega_{{i}}\}\in
M_{\omega}(\mathrm{extr}\mathfrak{S}(\mathcal{H}\otimes\mathcal{K}))}\sum_{i}\pi_{i}f\circ\Theta(\omega_{i}),
\quad\omega\in\mathfrak{S}(\mathcal{H}\otimes\mathcal{K}),
\end{equation}
where $\Theta: \omega\mapsto\mathrm{Tr}_{\mathcal{K}}\omega$ is the
partial trace \cite{H-QI} (by the spectral theorem the right hand
side of (\ref{E-f}) is well defined). If the von Neumann entropy
$H(\rho)=-\mathrm{Tr}\rho\log\rho$ is used as a function $f$, then
this method provides the construction of the Entanglement of
Formation, which is one of the most useful entanglement
measures\footnote{Entanglement measure is a EM having the particular
additional properties \cite{P&V}.} \cite{P&V}.

In what follows the properties of the function  $E^{f}$ defined by
(\ref{E-f}) in case of infinite dimensional spaces $\mathcal{H}$ and
$\mathcal{K}$ are considered.

An important problem in constructing of EM is to analyse the
continuity properties, in particular, to prove its continuity on the
entire state space $\mathfrak{S}(\mathcal{H}\otimes\mathcal{K})$
(formally the last property is not included in the definition of EM,
but in finite dimensions it is considered as a natural requirement).
Note that the continuity of the function $E^{f}$ is not obvious even
in the finite dimensional case and in general it is proved by using
the explicit form of the function $f$. By Theorem~\ref{main} the
$\mu$-compactness and stability of the set
$\mathfrak{S}(\mathcal{H}\otimes\mathcal{K})$ guarantee continuity
of the function $E^{f}$ on this set for an arbitrary continuous
function $f$ in the both finite and infinite dimensional
cases.\vspace{5pt}

\begin{theorem}\label{em}
\textit{Let $f$ be a concave continuous nonnegative function on the
set $\,\mathfrak{S}(\mathcal{H})$ satisfying condition
(\ref{f-property}). Then the function $E^{f}$ defined by (\ref{E-f})
is an entanglement monotone that is continuous on the set
$\,\mathfrak{S}(\mathcal{H}\otimes\mathcal{K})$.}
\end{theorem}
\vspace{5pt} \noindent {\tt Proof.} Nonnegativity, concavity, and
continuity of the function $f$ imply its boundedness. By
Theorem~\ref{main} the stability property of the $\mu$-compact set
 $\mathfrak{S}(\mathcal{H}\otimes\mathcal{K})$
guarantees the continuity of the function
$\mathrm{co}(f\circ\Theta)$ and therefore, its coincidence with
the function $\overline{\mathrm{co}}(f\circ\Theta)$, which by
Proposition~\ref{formula} has the following representation
\begin{equation}\label{co-cl-exp-em}
\overline{\mathrm{co}}(f\circ\Theta)(\omega) = \ \inf_{\mu\in
M_{\omega}(\mathfrak{S}(\mathcal{H}\otimes\mathcal{K}))}\int_{\mathfrak{S}(\mathcal{H}\otimes\mathcal{K})}
(f\circ\Theta)(\varpi) \, \mu(d\varpi)\, , \qquad \omega\, \in \,
\mathfrak{S}(\mathcal{H}\otimes\mathcal{K}),
\end{equation}
where the infimum is achieved at some particular measure
$\mu_{\omega}$ in
$M_{\omega}(\mathfrak{S}(\mathcal{H}\otimes\mathcal{K}))$. By
concavity, continuity, and boundedness of the function
$f\circ\Theta$ one can assume that $\mu_{\omega}$ is a measure in
$M_{\omega}(\mathrm{extr}\mathfrak{S}(\mathcal{H}\otimes\mathcal{K}))$.
Hence the definition of the function $E^{f}$ and concavity of the
function  $f\circ\Theta$ imply
$E^{f}=\mathrm{co}(f\circ\Theta)=\overline{\mathrm{co}}(f\circ\Theta)$.

By (\ref{f-property}) the nonnegative function $f\circ\Theta$
vanishes  on a pure state in
$\mathfrak{S}(\mathcal{H}\otimes\mathcal{K})$ if and only if this
state is product. As shown in  \cite{H-Sh-W}, a state $\omega$ is
nonentangled if and only if there exists such measure $\mu_{\omega}$
supported by pure product states that
$\mathbf{b}(\mu_{\omega})=\omega$. Thus the above remark shows that
condition  E-1 is fulfilled for the function
$E^{f}=\overline{\mathrm{co}}(f\circ\Theta)$.

Condition E-2 for the function $E^{f}$ is established in the same
way as in the finite dimensional case (see \cite{P&V}).

Condition E-3 for the function $E^{f}$ follows from its definition.

{\hfill $\Box$}
\medskip

\begin{example}\label{em+}
Generalizing the observation in \cite{Osborne} to the infinite
dimensional case, consider the family of continuous concave
functions
$$
f_{\alpha}(\rho)=2(1-\mathrm{Tr}\rho^{\alpha}),\quad \alpha>1,
$$
on the set $\mathfrak{S}(\mathcal{H})$ with
$\dim\mathcal{H}\leq+\infty$. It is easy to see that all the
functions in this family satisfy assumptions (\ref{f-property}). By
Theorem~\ref{em} $\{E^{f_{\alpha}}\}_{\alpha>1}$ is a family of
entanglement monotones, which are continuous and bounded on the set
$\mathfrak{S}(\mathcal{H}\otimes\mathcal{K})$ with
$\dim\mathcal{H}\leq+\infty$ and $\dim\mathcal{K}\leq+\infty$. The
case $\alpha=2$ is of special interest since the entanglement
monotone $E^{f_{2}}$ can be considered as an infinite dimensional
generalization of the notion of I-tangle \cite{R&C}.
\end{example}

\section{Possible generalizations and open questions}
\bigskip

Proposition \ref{ce-prop-0} actually gives an equivalent definition
of the $\mu$-compactness property for convex sets of the class
considered in this paper. A convex set $\mathcal{A}\, $ is
$\mu$-compact if and only if for arbitrary compact set
$\mathcal{K}\subseteq \mathcal{A}$ and for arbitrary $\varepsilon>0$
there exists a compact set $\mathcal{K}_{\varepsilon}\subseteq
\mathcal{A}$ such that for any expansion  of a point $\, x\in
\mathcal{K}$ into convex combination of points in $\, \mathcal{A}$
the total weight of points belonging to the set
$\mathcal{K}_{\varepsilon}$ is not less than $\, 1-\varepsilon$.
This property of (arbitrary!) convex set can be called generalized
$\, \mu$-compactness, or, in short, $\, \tilde
\mu$\nobreakdash-\hspace{0pt}compactness. By
Proposition~\ref{ce-prop-0} the $\, \tilde \mu$-compactness property
means $\, \mu$-compactness for convex bounded subsets of locally
convex spaces, that are complete separable metric spaces. The above
definition of $\, \tilde \mu$-compactness is translated without any
change to any convex closed subsets of linear topological spaces,
not necessarily bounded. The definition of $\, \mu$-compactness is,
in contrast,  not generalized to unbounded sets, since for unbounded
sets the barycenter map may  not be well defined: the integral $\,
\mathbf b (\mu ) \, = \, \int_\mathcal{A} \, x \, \mu(dx)\, $ may
not exist for some measures $\mu \in M(\mathcal{A})$. Thus the
notion of $\, \tilde \mu$-compactness generalizes the notion of $\,
\mu$-compactness to a wider class of convex sets. Similar to the
case of $\mu$-compact sets, the intersection and the Cartesian
product of finite or countable family of $\, \tilde \mu$-compact
sets is $\, \tilde \mu$-compact, the convex closed subset of a $\,
\tilde \mu$-compact set is $\, \tilde \mu$-compact. The proof of
this assertions is literally the same as the proof of
Proposition~\ref{conv-cl-mu-compact}. A complete analog of
Proposition~\ref{mu-compact} for continuous transformations of
$\mu$-compact sets holds as well. Nontrivial examples of $\, \tilde
\mu$-compacts sets appear even in the finite dimensional case.
\medskip

\begin{lemma}\label{varepsilon-d}
\textit{An arbitrary convex closed pointed (not containing any line)
cone in~$\mathbb  R^d$ is $\, \tilde \mu$-compact.}
\end{lemma}
\smallskip

\noindent {\tt Proof.} Let $\, \mathcal{C} \, \subset \, \mathbb
R^d$ be a convex pointed cone. Hence there exists a vector $\, a
\in \mathbb  R^d\, $ such that $\, \inf\limits_{x \in \mathcal C\,
, \|x\| \, = \, 1} \bigl( \, x \, , \, a\, \bigr) \, > \, 0$
\cite[p.53]{Boyd}. Then for each $r > 0$ the truncated cone
$$
\mathcal C_r \, = \, \bigl\{\, x \in \mathcal C \, , \, (x, a) \,
\le \, r \bigr\}
$$
is compact. Arbitrary compact set $\, \mathcal K \, \subset \,
\mathcal C\, $ can be embedded into a particular truncated cone $\,
\mathcal C_r$. Then for each $\, \varepsilon
> 0\, $ the compact set  $\, \mathcal K_{\varepsilon}\, = \, \mathcal C_{\,
\frac{r}{\varepsilon}}$ has the required property.

{\hfill $\Box$}
\medskip

The following result from the convex geometry is well known, so we
omit its proof.
\medskip

\begin{lemma}\label{three}
\textit{The following properties of a convex closed set $\,
\mathcal{A} \subset \mathbb R^d$ are equivalent:}
\begin{enumerate}[(i)]
 \item \textit{$\, \mathcal{A}\, $ is contained in a convex pointed cone;}
 \item \textit{$\, \mathcal{A}\, $ has at least one extreme point;}
 \item \textit{$\, \mathcal{A}\, $ does not contain a straight line;}
 \item \textit{the polar of the set $\, \mathcal{A}\, $ has a nonempty interior.}
\end{enumerate}
\end{lemma}
\medskip

By applying Lemma~\ref{varepsilon-d} we conclude that
property~$\mathrm{(i)}$ implies $\, \tilde \mu$-compactness of the
set $\mathcal{A}$. On the other hand, a line is not $\, \tilde
\mu$-compact, hence $\, \tilde \mu$-compactness implies
property~$\mathrm{(iii)}$. Thus, we obtain the following result.
\medskip

\begin{property}\label{final}
\textit{In the space $\,\mathbb R^d$ the $\, \tilde \mu$-compactness
is equivalent to each of properties $\mathrm{(i)} - \mathrm{(iv)}$
in Lemma~\ref{three}}.
\end{property}
\medskip

Positive cones in the spaces $l_{\, p}\, $ and $\, L_{\, p}(X)$ are
$\, \tilde \mu$-compact in the weak topology. The proof is analogous
to the proof of Lemma~\ref{varepsilon-d} and is based on the weak
compactness of bounded sets in this spaces. The positive cone in
$l_{\, 1}$ (Proposition~\ref{V1}), the cone of finite Borel measures
on a complete separable metric space (Proposition~\ref{measures}),
and the cone of positive operators in the Shatten class of order
$p=1$ (Proposition~\ref{Lp}) are $\tilde \mu$-compact and, in
contrast to these propositions, one need not take bounded parts of
these cones. Thus,  $\tilde \mu$\nobreakdash-\hspace{0pt}compactness
substantially extends the notion of $\mu$-compactness. This
motivates our  first question.
\medskip

\noindent \textbf{Question 1.} \textit{To what extent are the
results of this paper generalized to $\tilde \mu$-compact sets~?}
\medskip

\noindent The following questions concern $\, \mu$-compact sets.
\medskip

\noindent \textbf{Question 2.} \textit{Do there exist
$\mu$-compact noncompact sets in the spaces  $L_p$ and $l_p$ with
$p > 1$~?}
\medskip

\noindent \textbf{Question 3.} \textit{Consider a Banach lattice.
Under what conditions is the bounded part of the positive cone in it
$\, \mu$-compact~?}
\medskip

\noindent \textbf{Question 4.} \textit{Under what conditions on a
convex set $\mathcal A$, is the convex hull of an arbitrary
continuous bounded function on this set continuous~?}
\medskip

The latter property holds for stable  $\mu$-compact sets
(Theorem~\ref{main}), but it does not hold for stable pointwise
$\mu$-compact sets that are not $\mu$-compact
(Proposition~\ref{formula-}). The unit ball in the space $l_2$
possesses  this property (Remark~\ref{r3}), but the positive part
of this ball does not (Example~\ref{conv-cl-ex}).

\section{Appendix}

\subsection{The compactness criterion for subsets of the cone
$\mathfrak{L}_{+}(\mathcal{H},\mathcal{H}')$}

Let $\mathfrak{L}_{+}(\mathcal{H},\mathcal{H}')$ be the cone of
linear continuous positive maps from  the Banach space
$\mathfrak{T}(\mathcal{H})$ of trace-class operators in a separable
Hilbert space $\mathcal{H}$ into the similar Banach space
$\mathfrak{T}(\mathcal{H'})$. The compactness criterion for subsets
of this cone in the strong operator topology is presented in the
following lemma.\medskip
\begin{lemma}\label{comp-c}
\textup{1)} \textit{A closed bounded subset
$\mathfrak{L}_{0}\subseteq\mathfrak{L}_{+}(\mathcal{H},\mathcal{H}')$
is compact in the strong operator topology if in
$\mathfrak{S}(\mathcal{H})$ there exists such full rank state
$\sigma$ that $\{\Phi(\sigma)\}_{\Phi\in\mathfrak{L}_{0}}$ is a
compact subset of $\mathfrak{T}(\mathcal{H}')$. }

\textup{2)} \textit{If a subset
$\mathfrak{L}_{0}\subseteq\mathfrak{L}_{+}(\mathcal{H},\mathcal{H}')$
is compact in the strong operator topology then
$\{\Phi(\sigma)\}_{\Phi\in\mathfrak{L}_{0}}$ is a compact subset of
$\mathfrak{T}(\mathcal{H}')$ for any state $\sigma$ in
$\mathfrak{S}(\mathcal{H})$.}
\end{lemma}
\smallskip

\noindent {\tt Proof.} 1) Let $\{|i\rangle\}$ be the basis of
eigenvectors of the state $\sigma$ arranged in nonincreasing order
and $\mathcal{H}_{m}$ be the eigen subspace generated by the first
$m$ vectors of this basis.

Let $\{\Phi_{n}\}$ be an arbitrary sequence of maps in
$\mathfrak{L}_{0}$.

Show that for each $m$ for an arbitrary operator $A$ in
$\mathfrak{T}(\mathcal{H}_{m})$ there exists such subsequence
$\{\Phi_{n_{k}}\}$ that the sequence $\{\Phi_{n_{k}}(A)\}_{k}$
converges in $\mathfrak{T}(\mathcal{H}')$. Suppose first that $A\geq
0$. Since $A\in\mathfrak{T}(\mathcal{H}_{m})$ there exists such
$\lambda_{A}>0$ that $\lambda_{A}A\leq\sigma$. By the compactness
criterion for subsets of $\mathfrak{T}(\mathcal{H}')$ (see the
Appendix in \cite{Sh-H-3}) for arbitrary $\varepsilon>0$ there
exists finite rank projector
$P_{\varepsilon}\in\mathfrak{B}(\mathcal{H}')$ such that
$\mathrm{Tr}(I_{\mathcal{H}'}-P_{\varepsilon})\Phi(\sigma)<\varepsilon$,
and hence
$\mathrm{Tr}(I_{\mathcal{H}'}-P_{\varepsilon})\Phi(A)<\lambda_{A}^{-1}\varepsilon$
for all $\Phi\in\mathfrak{L}_{0}$. By the same compactness criterion
the set $\{\Phi(A)\}_{\Phi\in\mathfrak{L}_{0}}$ is compact. This
implies existence of the desired subsequence for a positive operator
$A$. Existence of such subsequence for an arbitrary operator
$A\in\mathfrak{T}(\mathcal{H}_{m})$ follows from representation of
this operator as a linear combination of positive operators in
$\mathfrak{T}(\mathcal{H}_{m})$.

Thus for each $m$ an arbitrary sequence
$\{\Phi_{n}\}\subset\mathfrak{L}_{0}$ contains such subsequence
$\{\Phi_{n_{k}}\}$ that there exists
\begin{equation}\label{lim-exp}
\lim_{k\rightarrow+\infty}\Phi_{n_{k}}(|i\rangle\langle
j|)=C^{m}_{ij}
\end{equation}
for all $i,j=\overline{1,m}$, where $\{C^{m}_{ij}\}$ are particular
operators in $\mathfrak{T}(\mathcal{H'})$.

For arbitrary  $m'>m$, by applying the above observation to the
sequence $\{\Phi_{n_{k}}\}_{k}$, we obtain such subsequence of the
sequence $\{\Phi_{n}\}$ that (\ref{lim-exp}) holds for all
$i,j=\overline{1,m'}$ with such set of operators  $\{C^{m'}_{ij}\}$,
that $C^{m'}_{ij}=C^{m}_{ij}$ for all $i,j=\overline{1,m}$.

By using this construction one can show existence of the set
$\{C_{ij}\}_{i,j=1}^{+\infty}$ of operators having the following
property: for each $m$ there exists such subsequence
$\{\Phi_{n_{k}}\}$ of the sequence $\{\Phi_{n}\}$ that
(\ref{lim-exp}) holds with $C^{m}_{ij}=C_{ij}$ for all
$i,j=\overline{1,m}$.

Consider the map $\Phi_{*}$ defined on the set
$\bigcup_{m\in\mathbb{N}}\mathfrak{T}(\mathcal{H}_{m})$  as follows
$$
\Phi_{*}\, :\, \sum_{i,j} a_{ij}\, |i\rangle\langle j|\quad \mapsto
\quad \sum_{i,j} a_{ij}\, C_{ij}\ \in\ \mathfrak{T}\,
(\mathcal{H'}).
$$
This map is linear by construction. It is easy to prove its
positivity and boundedness. Indeed, by the property of the set
$\{C_{ij}\}$ for an arbitrary operator
$A\in\bigcup_{m}\mathfrak{T}(\mathcal{H}_{m})$ there exists a
subsequence $\{\Phi_{n_{k}}\}$ of the sequence $\{\Phi_{n}\}$ such
that $\, \Phi_{*}(A)\, =\,
\lim\limits_{k\rightarrow+\infty}\Phi_{n_{k}}(A)$. Thus positivity
and boundedness of the map $\Phi_{*}$ follows from positivity of the
maps in the sequence $\{\Phi_{n}\}$ and uniform boundedness of these
maps. Since the set $\bigcup_{m}\mathfrak{T}(\mathcal{H}_{m})$ is
dense in $\mathfrak{T}(\mathcal{H})$, the map $\Phi_{*}$ can be
extended to a linear positive bounded map from
$\mathfrak{T}(\mathcal{H})$ into $\mathfrak{T}(\mathcal{H}')$
(denoted by the same symbol $\Phi_{*}$).

Show that the map $\Phi_{*}$ is a limit point of the sequence
$\{\Phi_{n}\}$ in the strong operator topology. This topology on
bounded subsets of $\mathfrak{L}_{+}(\mathcal{H},\mathcal{H}')$ can
be determined by countable family of seminorms
$\Phi\mapsto\|\Phi(\rho_{i})\|_{1}$, where $\{\rho_{i}\}$ is an
arbitrary countable dense subset of the set
$\mathfrak{S}(\mathcal{H})$.\footnote{Here the possibility to
express an arbitrary operator in $\mathfrak{T}(\mathcal{H})$ as a
linear combination of four states in $\mathfrak{S}(\mathcal{H})$ is
used.} We choose the set of states in
$\bigcup_{m}\mathfrak{T}(\mathcal{H}_{m})$ on the role of this
subset. An arbitrary vicinity of the map $\Phi_{*}$ contains
vicinity of the form
$$
\left\{\Phi\in\mathfrak{L}(\mathcal{H},\mathcal{H}')\,|\,\|(\Phi-\Phi_{*})(\rho_{i_{t}})\|_{1}<\varepsilon,
t=\overline{1,p}\right\},
$$
where $\{\rho_{i_{t}}\}_{t=1}^{p}$ is a finite subset of the above
set of states and $\varepsilon>0$. Since
$\{\rho_{i_{t}}\}_{t=1}^{p}\subset\mathfrak{T}(\mathcal{H}_{m})$ for
a particular $m$, the construction of the map $\Phi_{*}$ implies
existence of such subsequence $\{\Phi_{n_{k}}\}$ of the sequence
$\{\Phi_{n}\}$ that
$\Phi_{*}(\rho_{i_{t}})=\lim_{k\rightarrow+\infty}\Phi_{n_{k}}(\rho_{i_{t}})$
for all $t=\overline{1,p}$. Hence at least one element of the
sequence $\{\Phi_{n}\}$ is contained in the above vicinity.

Thus the map $\Phi_{*}$ is a limit point of the sequence
$\{\Phi_{n}\}$ in the strong operator topology. By metrizability of
the strong operator topology on bounded subsets of the cone
$\mathfrak{L}_{+}(\mathcal{H},\mathcal{H}')$ this implies existence
of a subsequence of the sequence $\{\Phi_{n}\}$ converging to the
map $\Phi_{*}$. Compactness of the set $\mathfrak{L}_{0}$ is proved.

2) This assertion immediately follows from the definition of the
strong operator topology.\medskip

{\hfill $\Box$} \medskip

\subsection{The openness criterion}

\begin{lemma}\label{open-map}
\textit{Let $\varphi$ be a map from a topological space $X$ to a
metric space $Y$. The following assertions are equivalent:}
\begin{enumerate}[(i)]
  \item \textit{the map $\varphi$ is open;}
  \item \textit{for arbitrary $x_{0}\in X$ and arbitrary sequence
$\{y_{n}\}\subset Y$, converging to $y_{0}=\varphi(x_{0})$, there
exists a subnet $\{y_{n_{\lambda}}\}_{\lambda\in\Lambda}$ of the
sequence $\{y_{n}\}$ and a net
$\{x_{\lambda}\}_{\lambda\in\Lambda}$, converging to $x_{0}$, such
that $\varphi(x_{\lambda})=y_{n_{\lambda}}$ for all
$\lambda\in\Lambda$.}
\end{enumerate}
\end{lemma}
\smallskip

\noindent {\tt Proof.} $\mathrm{(i)}\Rightarrow\mathrm{(ii)}$ Let
$\mathfrak{U}$ be the set of all vicinities of the point $x_{0}$.
Then the set  $\Lambda$ of all pairs $\lambda=(U,k)$, where
$U\in\mathfrak{U}$ and $k\in\mathbb{N}$, with the partial order
$$
\{\lambda_{1}=(U_{1},k_{1})\succ\lambda_{2}=(U_{2},k_{2})\}\;\Leftrightarrow\;
\{k_{1}\geq k_{2}\;\; \textrm{and}\;\; U_{1}\subseteq U_{2}\}
$$
is directed. For each $\lambda=(U,k)$ the set
$W_{\lambda}=\varphi(U)\cap V_{k}$, where $V_{k}$ is the open ball
in $Y$ with the center $y_{0}$ and radius  $1/k$, is a vicinity of
the point $y_{0}$. Hence there exists minimal natural $n_{\lambda}$
such that $y_{n_{\lambda}}\in W_{\lambda}$. It is easy to see that
$\{y_{n_{\lambda}}\}_{\lambda\in\Lambda}$ is a subnet of the
sequence $\{y_{n}\}$. For each $\lambda=(U,k)$ there exists such
$x_{\lambda}\in U$ that $\varphi(x_{\lambda})=y_{n_{\lambda}}$. It
is clear that the net $\{x_{\lambda}\}_{\lambda\in\Lambda}$
converges to $x_{0}$.

$\mathrm{(ii)}\Rightarrow\mathrm{(i)}$ If there exists such open set
$U\subseteq X$ that the set  $\varphi(U)$ is not open then there
exist $y_{0}=\varphi(x_{0})\in\varphi(U)$ and sequence
$\{y_{n}\}\subset Y\setminus\varphi(U)$ converging to $y_{0}$. By
using $\mathrm{(ii)}$ it is easy to obtain a contradiction.

{\hfill $\Box$} \bigskip

The authors are grateful to the referees for useful remarks and
recommendations providing improvement of this paper. \medskip

The research of the first author is supported by the grants RFBR
08-01-00208,\break MD-2195.2008.1, and NSh-3233.2008.1, the second
author is supported by the grants RFBR 07-01-00156 and 09-01-00424a
and by the analytical departmental target program "Development of
scientific potential of the higher school 2009-2010", grant
2.1.1/500.

\end{document}